\documentclass{article}%
\usepackage{amsfonts}
\usepackage{amsmath}
\usepackage{amssymb}
\usepackage{amsxtra}
\usepackage{geometry}
\usepackage{graphicx}%
\providecommand{\U}[1]{\protect\rule{.1in}{.1in}}
\newtheorem{theorem}{Theorem}

\newtheorem{corollary}[theorem]{Corollary}

\newtheorem{lemma}[theorem]{Lemma}

\newtheorem{proposition}[theorem]{Proposition}

\numberwithin{equation}{section}
\begin{document}

\title{Spectral approach to homogenization of an elliptic operator periodic in some directions.}
\author{R. Bunoiu\\LMAM, UMR 7122, Universit\'{e} de Metz et CNRS\\Ile du Saulcy, F-57045 METZ Cedex 1, France\\email: bunoiu@math.univ-metz.fr
\and G. Cardone\\University of Sannio - Department of Engineering\\Piazza Roma, 21 - 84100 Benevento, Italy \\email: giuseppe.cardone@unisannio.it
\and T. Suslina\\St. Petersburg State University, Department of Physics \\Petrodvorets, Ul'yanovskaya 1, St. Petersburg 198504, Russia\\email: suslina@list.ru}
\maketitle

\begin{abstract}
The operator
\[
A_{\varepsilon}= D_{1} g_{1}(x_{1}/\varepsilon, x_{2}) D_{1} + D_{2}
g_{2}(x_{1}/\varepsilon, x_{2}) D_{2}
\]
is considered in $L_{2}({\mathbb{R}}^{2})$, where $g_{j}(x_{1},x_{2})$,
$j=1,2,$ are periodic in $x_{1}$ with period 1, bounded and positive definite.
Let function $Q(x_{1},x_{2})$ be bounded, positive definite and periodic in
$x_{1}$ with period 1. Let $Q^{\varepsilon}(x_{1},x_{2})= Q(x_{1}/\varepsilon,
x_{2})$. The behavior of the operator $(A_{\varepsilon}+ Q^{\varepsilon}%
)^{-1}$ as $\varepsilon\to0$ is studied. It is proved that the operator
$(A_{\varepsilon}+ Q^{\varepsilon})^{-1}$ tends to $(A^{0} + Q^{0})^{-1}$ in
the operator norm in $L_{2}(\mathbb{R}^{2})$. Here $A^{0}$ is the effective
operator whose coefficients depend only on $x_{2}$, $Q^{0}$ is the mean value
of $Q$ in $x_{1}$. A sharp order estimate for the norm of the difference
$(A_{\varepsilon}+ Q^{\varepsilon})^{-1}- (A^{0} + Q^{0})^{-1}$ is obtained.
The result is applied to homogenization of the Schr\"odinger operator with a
singular potential periodic in one direction.

\end{abstract}

\section{Introduction\label{sect0}}

\subsection{Spectral approach to homogenization problems}

Homogenization problems for periodic differential operators of mathematical
physics are of significant interest both from theoretical point of view and
for applications. A broad literature is devoted to homogenization problems. At
the first place, the books \cite{BeLP}, \cite{Sa}, \cite{BaPa}, \cite{ZhKO}
should be mentioned.

Consider a typical homogenization problem. Let $\mathcal{A}_{\varepsilon}$ be
a family of operators in $L_{2}(\mathbb{R}^{d})$ given by
\begin{equation}
{\mathcal{A}}_{\varepsilon} = - \hbox{div}\,g \left(  {\mathbf{x}%
}/{\varepsilon}\right)  \nabla, \quad\varepsilon>0. \label{0.1}%
\end{equation}
Here $g(\mathbf{x})$ is a positive definite and bounded $(d \times
d)$-matrix-valued function. It is assumed that $g(\mathbf{x})$ is periodic
with respect to some lattice $\Gamma\subset\mathbb{R}^{d}$. The problem is to
study the behavior of the solution $u_{\varepsilon}$ of the equation
\begin{equation}
\mathcal{A}_{\varepsilon}u_{\varepsilon}+ u_{\varepsilon}= F,\quad F \in
L_{2}(\mathbb{R}^{d}), \label{0.2}%
\end{equation}
for small $\varepsilon$. It turns out that $u_{\varepsilon}$ tends (in some
sense) to $u_{0}$ as $\varepsilon\to0$, where $u_{0}$ is the solution of the
"homogenized" equation
\begin{equation}
\mathcal{A}^{0} u_{0} + u_{0} = F. \label{0.3}%
\end{equation}
Here $\mathcal{A}^{0} = -\hbox{div}\, g^{0} \nabla$, and $g^{0}$ is a constant
positive matrix called the \textit{effective matrix}. The operator
$\mathcal{A}^{0}$ is called the \textit{effective operator}. The character of
convergence of $u_{\varepsilon}$ to $u_{0}$ and error estimates are studied.

A possible approach to homogenization of periodic differential operators in
$\mathbb{R}^{d}$ consists in application of the Floquet-Bloch theory combining
with methods of the analytic perturbation theory. In this connection, we
mention the papers \cite{Se}, \cite{Zh}, \cite{CoV}. In a series of papers
\cite{BSu1}, \cite{BSu2}, \cite{BSu3}, \cite{BSu4}, a new operator-theoretic
approach to homogenization problems was suggested and developed. It was also
based on the Floquet-Bloch decomposition for a periodic operator and applying
methods of the analytic perturbation theory. The main idea of this approach
was to study a homogenization procedure as a spectral \textit{threshold
effect} at the bottom of the spectrum of a periodic elliptic operator. In
particular, this approach allowed the authors to prove the following error
estimate:
\begin{equation}
\| u_{\varepsilon}-u_{0}\|_{L_{2}(\mathbb{R}^{d})} \le C\varepsilon
\|F\|_{L_{2}(\mathbb{R}^{d})}, \label{0.4}%
\end{equation}
see \cite{BSu1}. This estimate is order-sharp and the constant $C$ is well
controlled. Estimate (\ref{0.4}) means that the resolvent $(\mathcal{A}%
_{\varepsilon}+I)^{-1}$ tends to $(\mathcal{A}^{0} + I)^{-1}$ in the operator
$L_{2}$-norm, and
\begin{equation}
\|(\mathcal{A}_{\varepsilon}+I)^{-1} -(\mathcal{A}^{0} +I)^{-1} \|_{L_{2}%
(\mathbb{R}^{d}) \to L_{2}(\mathbb{R}^{d})} \le C \varepsilon. \label{0.5}%
\end{equation}
Such estimates are called "operator error estimates"; they were obtained in
\cite{BSu1} for a wide class of matrix periodic differential operators.

Let us briefly discuss the method of the proof of estimate (\ref{0.5}). By the
scaling transformation, (\ref{0.5}) is equivalent to the estimate
\begin{equation}
\|(\mathcal{A} +\varepsilon^{2} I)^{-1} -(\mathcal{A}^{0} +\varepsilon^{2}
I)^{-1} \|_{L_{2}(\mathbb{R}^{d}) \to L_{2}(\mathbb{R}^{d})} \le C
\varepsilon^{-1}. \label{0.6}%
\end{equation}
Here $\mathcal{A} = -\hbox{div}\,g(\mathbf{x}) \nabla$. The bottom of the
spectrum of $\mathcal{A}$ is the point $\lambda=0$. It is natural that the
behavior of the resolvent $(\mathcal{A} +\varepsilon^{2} I)^{-1}$ can be
described in terms of the threshold characteristics of $\mathcal{A}$ (i.~e.,
the spectral characteristics at the bottom of the spectrum). Applying the
Floquet-Bloch theory, we decompose $\mathcal{A}$ in the direct integral of
operators $\mathcal{A}(\mathbf{k})$ acting in $L_{2}(\Omega)$. Here $\Omega$
is the cell of the lattice $\Gamma$. The parameter $\mathbf{k}\in
\mathbb{R}^{d}$ is called a \textit{quasimomentum}. The operator
$\mathcal{A}(\mathbf{k})$ is given by the differential expression $(\mathbf{D}
+ \mathbf{k})^{*}g(\mathbf{x})(\mathbf{D} + \mathbf{k})$ with periodic
boundary conditions. Estimate (\ref{0.6}) is equivalent to the estimate
\begin{equation}
\|(\mathcal{A}(\mathbf{k}) +\varepsilon^{2} I)^{-1} -(\mathcal{A}%
^{0}(\mathbf{k}) +\varepsilon^{2} I)^{-1} \|_{L_{2}(\Omega) \to L_{2}(\Omega)}
\le C \varepsilon^{-1}, \label{0.7}%
\end{equation}
which must be uniform in $\mathbf{k} \in\widetilde{\Omega}$. Here
$\widetilde{\Omega}$ is the Brillouin zone of the dual lattice $\widetilde
{\Gamma}$. The main part of investigation is the study of the operator family
$\mathcal{A}(\mathbf{k})$ by means of the analytic perturbation theory. A
crucial role is played by the \textit{spectral germ} of the operator family
$\mathcal{A}(\mathbf{k})$ at $\mathbf{k}=0$. The spectral germ is a finite
rank selfadjoint operator which is defined in terms of the threshold
characteristics of $\mathcal{A}$. It is possible to find a finite rank
approximation for the resolvent $(\mathcal{A}(\mathbf{k}) +\varepsilon^{2}
I)^{-1}$ for small $\varepsilon$ in terms of the spectral germ. Next, it turns
out that the family $\mathcal{A}^{0}(\mathbf{k})$ corresponding to the
effective operator \textit{has the same spectral germ} as the family
$\mathcal{A}(\mathbf{k})$. This allows one to pass from approximation in terms
of the germ to approximation in terms of the effective operator.

\subsection{Homogenization problems for operators periodic not in all
directions}

Along with operators whose coefficients are periodic in all directions, it is
interesting to study operators with coefficients that are periodic only in
part of variables. Such problems are important for applications, in
particular, in the theory of waveguides. In these problems, the coefficients
of the effective operator depend on the "nonperiodic" variables, while
homogenization is related to the "periodic" variables. Technically, the study
of such problems is more difficult than the study of operators with
coefficients that are periodic in all directions.

In the paper \cite{Su}, the following elliptic operator in a strip
$\Pi=\mathbb{R}\times(0,a)$ was considered:
\[
{\mathcal{B}}_{\varepsilon}=D_{1}g_{1}(\varepsilon^{-1}x_{1},x_{2})D_{1}%
+D_{2}g_{2}(\varepsilon^{-1}x_{1},x_{2})D_{2}.
\]
The coefficients $g_{1},g_{2}$ were assumed to be positive definite, bounded
and 1-periodic in $x_{1}$. Moreover, it was assumed that $g_{1},g_{2}$ are
periodic Lipschitz functions with respect to $x_{2}$. On the boundary
$\partial\Pi$, periodic boundary conditions were posed. It was shown that
there exists an "effective operator" ${\mathcal{B}}^{0}=D_{1}g_{1}^{0}%
(x_{2})D_{1}+D_{2}g_{2}^{0}(x_{2})D_{2}$ (also with periodic boundary
conditions on $\partial\Pi$) such that
\begin{equation}
\Vert(\mathcal{B}_{\varepsilon}+I)^{-1}-(\mathcal{B}^{0}+I)^{-1}\Vert
_{L_{2}(\Pi)\rightarrow L_{2}(\Pi)}\leq C\varepsilon. \label{0.8}%
\end{equation}
As before, the scaling transformation is used and the Floquet decomposition
(in $x_{1}$) is applied. The quasimomentum $k$ is one-dimensional. Then the
problem reduces to the study of the operator family
\[
{\mathcal{B}}(k,\varepsilon)=\mathcal{B}_{1}(k)+\varepsilon^{2}\mathcal{B}%
_{2},
\]
acting in $L_{2}(\Omega)$, where $\Omega=(0,1)\times(0,a)$. Here
${\mathcal{B}}_{1}(k)$ is given by the expression $(D_{1}+k)g_{1}%
(\mathbf{x})(D_{1}+k)$ and $\mathcal{B}_{2}=D_{2}g_{2}(\mathbf{x})D_{2}$; the
boundary conditions are periodic in both variables. The method of \cite{BSu1}
can be applied to the operator $\mathcal{B}_{1}(k)$ "layerwise". This allows
to introduce the spectral germ for the family $\mathcal{B}_{1}(k)$ and to
obtain approximation for $\mathcal{B}_{1}(k)$ in terms of the germ. Main
technical difficulties are related to adding the \textit{unbounded}
"perturbation" term $\varepsilon\mathcal{B}_{2}$ to $\mathcal{B}_{1}(k)$ and
approximating the operator $(\mathcal{B}(k,\varepsilon)+\varepsilon^{2}%
I)^{-1}$.

\subsection{Main goal and main results}

In the present paper, we study an elliptic operator
\begin{equation}
{A}_{\varepsilon}= D_{1} g_{1}(\varepsilon^{-1}x_{1},x_{2})D_{1} + D_{2}
g_{2}(\varepsilon^{-1}x_{1},x_{2})D_{2}, \quad\varepsilon>0, \label{0.9}%
\end{equation}
acting in $L_{2}(\mathbb{R}^{2})$. The coefficients $g_{j}(x_{1},x_{2})$,
$j=1,2,$ are assumed to be bounded, Lipschitz in $x_{2}$, positive definite
and $1$-periodic in $x_{1}$. Next, let $Q(\mathbf{x})$ be a positive definite
and bounded function in $\mathbb{R}^{2}$. It is assumed that $Q$ is 1-periodic
in $x_{1}$ and Lipschitz in $x_{2}$. We denote $Q^{\varepsilon}(\mathbf{x}) =
Q(\varepsilon^{-1}{x}_{1}, x_{2})$.

\textit{Our main goal} is to find an approximation in the operator norm in
$L_{2}(\mathbb{R}^{2})$ for the \textit{generalized resolvent}
$(A_{\varepsilon}+ Q^{\varepsilon})^{-1}$ of the operator (\ref{0.9}).

We construct the effective operator
\[
A^{0} = D_{1} g_{1}^{0}(x_{2}) D_{1} + D_{2} g_{2}^{0}(x_{2}) D_{2}
\]
with the effective coefficients depending only on $x_{2}$ (see (\ref{1.11}),
(\ref{1.12})). By $Q^{0}(x_{2})$ we denote the mean value of $Q(x_{1},x_{2})$
in $x_{1}$.

The \textit{main result} of the paper is the following estimate
\begin{equation}
\|({A}_{\varepsilon}+Q^{\varepsilon})^{-1} - ({A}^{0}+Q^{0})^{-1}%
\|_{L_{2}(\mathbb{R}^{2}) \to L_{2}(\mathbb{R}^{2})} \le C \varepsilon.
\label{0.10}%
\end{equation}
By the scaling transformation, the question is reduced to the study of the
operator family $(A(\varepsilon) + \varepsilon^{2} Q)^{-1}$, where the
operator
\begin{equation}
A(\varepsilon) = D_{1}g_{1}(x_{1},x_{2})D_{1} + \varepsilon^{2} D_{2}
g_{2}(x_{1},x_{2})D_{2} \label{0.11}%
\end{equation}
contains small factor $\varepsilon^{2}$ standing at the second order operator
$D_{2} g_{2} D_{2}$. This is the main difficulty. We apply the Floquet-Bloch
decomposition (in $x_{1}$) and decompose the operator (\ref{0.11}) in the
direct integral of operators
\begin{equation}
A(k,\varepsilon) = A_{1}(k) + \varepsilon^{2} A_{2}, \label{0.12}%
\end{equation}
acting in $L_{2}(\Omega)$. Now the domain $\Omega$ is \textit{unbounded}:
$\Omega= (0,1)\times\mathbb{R}$. The operator $A_{1}(k)$ is given by the
expression $(D_{1}+k) g_{1}(\mathbf{x}) (D_{1}+k)$ with periodic boundary
conditions on $\partial\Omega$, and $A_{2} = D_{2} g_{2}(\mathbf{x})D_{2}$. We
apply method of \cite{BSu1} to the operator $A_{1}(k)$ "layerwise" (for each
fixed $x_{2}$), and define the spectral germ of the family $A_{1}(k)$ at
$k=0$. Now the germ is not a finite rank operator. It is possible to
approximate the part of $A_{1}(k)$ near the bottom of the spectrum in terms of
the germ. Main difficulties are related to taking the term $\varepsilon^{2}
A_{2}$ into account.

The essential differences of the problem studied here from the problem studied
before in \cite{Su} is that $\Omega$ is \textit{unbounded} domain and that we
study the \textit{generalized} resolvent of $A_{\varepsilon}$ (instead of the
ordinary resolvent). We develop some technical tools suggested in \cite{Su}
and adapt them for the case of unbounded domain $\Omega$ and generalized resolvent.

Main results can be generalized for the case of arbitrary dimension without
any new technical difficulties. We give formulations of the corresponding
results (see Section \ref{sect11}). In the main exposition, we consider the
case $d=2$ only for the sake of simplicity and clearness.

Our general result can be applied to homogenization of the Schr\"{o}dinger
operator with a singular potential periodic in one direction and non-periodic
in another direction; see Section \ref{sect10}. This operator arises in the
model of a "soft waveguide" (the spectral properties of such operators were
studied in \cite{FK}). We emphasize that for application to the
Schr\"{o}dinger operator we need to study in advance a generalized resolvent
of $A_{\varepsilon}$ (but not the ordinary resolvent).

\subsection{Plan of exposition}

In Section \ref{sect1}, we give the precise statement of the problem and
formulate the main result (Theorem 1). In Section \ref{sect2}, using the
scaling transformation and the Floquet-Bloch decomposition we reduce the
problem to the study of the generalized resolvent $(A(k,\varepsilon
)+\varepsilon^{2}Q)^{-1}$ of the operator (\ref{0.12}). In Section
\ref{sect3}, the operator family $A_{1}(k)$ is studied by using method of
\cite{BSu1}; the spectral germ of $A_{1}(k)$ is introduced. In Section
\ref{sect4}, Theorem 1 is reduced to Theorem \ref{theor4.1}, which gives
approximation of the operator $(A(k,\varepsilon)+\varepsilon^{2}Q)^{-1}$ in
terms of the spectral germ. Theorem \ref{theor4.1} is proved in Sections
\ref{sect5}--\ref{sect7}; this is technical part of the paper. Section
\ref{sect10} is devoted to application of Theorem \ref{theor1.1} to
homogenization problem for the Schr\"{o}dinger operator with a singular
potential periodic in one direction. Finally, in Section \ref{sect11}
generalizations of main results for the case of arbitrary dimension are formulated.

\subsection{Notation}

Let $\mathfrak{H}$ and $\mathfrak{G}$ be separable Hilbert spaces. The symbols
$(\cdot,\cdot)_{\mathfrak{H}}$, $\| \cdot\|_{\mathfrak{H}}$ stand for the
inner product and the norm in $\mathfrak{H}$; the symbol $\|T\|_{\mathfrak{H}%
\to\mathfrak{G}}$ denotes the norm of a bounded operator $T: \mathfrak{H}
\to\mathfrak{G}$. Sometimes we omit the indices if this does not lead to
confusion. If $\mathfrak{N}$ is a subspace in $\mathfrak{H}$, then its
orthogonal complement is denoted by $\mathfrak{N}^{\perp}$. If $P$ is the
orthogonal projection of $\mathfrak{H}$ onto $\mathfrak{N}$, then $P^{\perp}$
denotes the orthogonal projection of $\mathfrak{H}$ onto $\mathfrak{N}^{\perp
}$.

By $\mathbf{1}_{d}$ we denote the unit $(d\times d)$-matrix. The symbol
$\langle\cdot, \cdot\rangle_{\mathbb{C}^{d}}$ stands for the standard inner
product in $\mathbb{C}^{d}$.

By $H^{s}(\mathcal{D})$, $s \in\mathbb{R}$, we denote the Sobolev spaces in a
domain $\mathcal{D} \subset\mathbb{R}^{d}$.

\section{Definition of the operators. Main results\label{sect1}}

\subsection{Statement of the problem}

Let $\mathbf{x}=(x_{1},x_{2})\in{\mathbb{R}}^{2}$. We use the notation
$\partial_{j}=\frac{\partial}{\partial x_{j}},\,D_{j}=-i\partial_{j},\,j=1,2.$
Assume that $g_{j}(\mathbf{x}),$ $j=1,2,$ are real-valued measurable functions
in ${\mathbb{R}}^{2}$ such that%
\begin{equation}
0<c_{0}\leq g_{j}(\mathbf{x})\leq c_{1}<\infty,\quad j=1,2,\quad
\hbox {a. e.}\ \mathbf{x}\in{\mathbb{R}}^{2}. \label{1.1}%
\end{equation}
The functions $g_{j}(x_{1},x_{2})$ are assumed to be periodic in $x_{1}$ with
period $1$:%
\begin{equation}
g_{j}(x_{1}+1,x_{2})=g_{j}(x_{1},x_{2}),\quad\mathbf{x}\in{\mathbb{R}}%
^{2},\ j=1,2. \label{1.2}%
\end{equation}
Suppose also that $g_{j}(x_{1},x_{2})$ are Lipschitz class with respect to
$x_{2}$:%
\begin{equation}
\underset{\mathbf{x}\in{\mathbb{R}}^{2}}{\operatorname*{ess}\sup}|\partial
_{2}g_{j}(\mathbf{x})|\leq c_{2}<\infty,\quad j=1,2. \label{1.3}%
\end{equation}
In the space $L_{2}({\mathbb{R}}^{2})$ we consider the operator
$A_{\varepsilon}$ given formally by the differential expression (\ref{0.9}).
Precisely, $A_{\varepsilon}$ is defined as a selfadjoint operator in
$L_{2}({\mathbb{R}}^{2})$ generated by the closed quadratic form%
\begin{equation}
a_{\varepsilon}[u,u]=\int_{{\mathbb{R}}^{2}}\left(  g_{1}\left(  \frac{x_{1}%
}{\varepsilon},x_{2}\right)  |D_{1}u|^{2}+g_{2}\left(  \frac{x_{1}%
}{\varepsilon},x_{2}\right)  |D_{2}u|^{2}\right)  d\mathbf{x},\text{ \ }u\in
H^{1}({\mathbb{R}}^{2}),\text{\ }\varepsilon>0. \label{1.6}%
\end{equation}
Let $Q(\mathbf{x})=Q(x_{1},x_{2})$ be a real-valued measurable function in
${\mathbb{R}}^{2}$ such that%
\begin{align}
Q(x_{1}+1,x_{2})  &  =Q(x_{1},x_{2}),\quad\mathbf{x}\in{\mathbb{R}}%
^{2},\label{1.6a}\\
0<c_{3}\leq Q(\mathbf{x})  &  \leq c_{4}<\infty,\quad\hbox {a. e.}\ \mathbf{x}%
\in{\mathbb{R}}^{2}. \label{1.7}%
\end{align}
Suppose also that $Q(\mathbf{x})$ is Lipschitz class with respect to $x_{2}$:%
\begin{equation}
\underset{\mathbf{x}\in{\mathbb{R}}^{2}}{\operatorname*{ess}\sup}|\partial
_{2}Q(\mathbf{x})|\leq c_{5}<+\infty. \label{1.8}%
\end{equation}
We denote
\begin{equation}
Q^{\varepsilon}(\mathbf{x})=Q\left(  \frac{x_{1}}{\varepsilon},x_{2}\right)  .
\label{1.9}%
\end{equation}

\textit{Our goal is to study the behavior of the generalized resolvent
$(A_{\varepsilon}+Q^{\varepsilon})^{-1}$ as $\varepsilon\rightarrow0.$}

\subsection{Main result}

We find the \textit{effective operator} $A^{0}$ of the same form as $A$, but
with coefficients independent of $x_{1}$, and the function $Q^{0}(x_{2})$ such
that
\[
(A_{\varepsilon}+Q^{\varepsilon})^{-1}\rightarrow(A^{0}+Q^{0})^{-1}%
\ \ \ \text{as }\varepsilon\rightarrow0,
\]
where convergence is understood in the operator norm in $L_{2}({\mathbb{R}%
}^{2})$.

The effective operator $A^{0}$ is a selfadjoint operator in $L_{2}%
({\mathbb{R}}^{2})$ generated by the quadratic form%
\begin{equation}
a^{0}[u,u]=\int_{{\mathbb{R}}^{2}}(g_{1}^{0}(x_{2})|D_{1}u|^{2}+g_{2}%
^{0}(x_{2})|D_{2}u|^{2})\,d\mathbf{x},\,\quad u\in H^{1}({\mathbb{R}}^{2}).
\label{1.10}%
\end{equation}
Here%
\begin{equation}
g_{1}^{0}(x_{2})=\left(  \int_{0}^{1}g_{1}(x_{1},x_{2})^{-1}dx_{1}\right)
^{-1}, \label{1.11}%
\end{equation}%
\begin{equation}
g_{2}^{0}(x_{2})=\int_{0}^{1}g_{2}(x_{1},x_{2})\,dx_{1}. \label{1.12}%
\end{equation}
Note that conditions (\ref{1.1}), (\ref{1.3}) imply the following inequalities
for the effective coefficients (\ref{1.11}), (\ref{1.12}):
\begin{equation}
0<c_{0}\leq g_{j}^{0}(x_{2})\leq c_{1}<+\infty,\quad x_{2}\in{\mathbb{R}},
\quad j=1,2, \label{1.13}%
\end{equation}%
\begin{align}
\underset{x_{2}\in{\mathbb{R}}}{\operatorname*{ess}\sup}|\partial_{2}g_{1}%
^{0}(x_{2})|  &  \leq c_{2}\left(  \frac{c_{1}}{c_{0}}\right)  ^{2}%
,\label{1.14}\\
\underset{x_{2}\in{\mathbb{R}}}{\operatorname*{ess}\sup}|\partial_{2}g_{2}%
^{0}(x_{2})|  &  \leq c_{2}. \label{1.14a}%
\end{align}
Due to conditions (\ref{1.14}), (\ref{1.14a}), $A^{0}$ can be given by the
differential expression $A^{0}=g_{1}^{0}(x_{2})D_{1}^{2}+D_{2}g_{2}^{0}%
(x_{2})D_{2}$ on the domain $H^{2}({\mathbb{R}}^{2}).$

Next, let $Q^{0}(x_{2})$ be the mean value of $Q(x_{1},x_{2})$:%
\begin{equation}
Q^{0}(x_{2})=\int_{0}^{1}Q(x_{1},x_{2})dx_{1}. \label{1.15}%
\end{equation}
From (\ref{1.7}) and (\ref{1.8}) it follows that%
\begin{align}
0  &  <c_{3}\leq Q^{0}(x_{2})\leq c_{4}<+\infty,\quad x_{2} \in{\mathbb{R}%
},\label{1.16}\\
&  \underset{x_{2}\in{\mathbb{R}}}{\operatorname*{ess}\sup}|\partial_{2}%
Q^{0}(x_{2})|\leq c_{5}<+\infty. \label{1.17}%
\end{align}

\textit{The main result} is given in the following theorem.

\begin{theorem}
\label{theor1.1}Suppose that conditions \textrm{(\ref{1.1})--(\ref{1.3})} are
satisfied and $A_{\varepsilon}$ is the operator in $L_{2}({\mathbb{R}}^{2})$
that corresponds to the form \textrm{(\ref{1.6})}. Let $A^{0}$ be the operator
corresponding to the form \textrm{(\ref{1.10})}, where the coefficients are
defined by \textrm{(\ref{1.11}), (\ref{1.12})}. Let $Q(\mathbf{x})$ be a
function in ${\mathbb{R}}^{2}$ satisfying \textrm{(\ref{1.6a})--(\ref{1.8})},
and $Q^{\varepsilon}$ is defined by \textrm{(\ref{1.9})}. Let $Q^{0}$ be
defined by \textrm{(\ref{1.15})}. Then we have%
\begin{equation}
\|(A_{\varepsilon}+Q^{\varepsilon})^{-1}-(A^{0}+Q^{0})^{-1}\|_{L_{2}%
({\mathbb{R}}^{2})\rightarrow L_{2}({\mathbb{R}}^{2})}\leq C\varepsilon
,\ \ 0<\varepsilon\leq1. \label{1.18}%
\end{equation}
The constant $C$ depends only on $c_{j},\,j=0,...,5.$
\end{theorem}

\section{Reduction to operators in the strip $\Omega$\label{sect2}}

\subsection{Scaling transformation}

We denote by $T_{\varepsilon},\,\varepsilon>0$, the unitary scaling
transformation in $L_{2}({\mathbb{R}}^{2})$ defined by the formula
\[
(T_{\varepsilon}u)(x_{1},x_{2})=\varepsilon^{\frac{1}{2}}u(\varepsilon
x_{1},x_{2}),\text{ \ }(x_{1},x_{2})\in{\mathbb{R}}^{2}.
\]
Let $A(\varepsilon)$ be the operator in $L_{2}({\mathbb{R}}^{2})$ generated by
the quadratic form%
\begin{equation}
a(\varepsilon)[u,u]=\int_{{\mathbb{R}}^{2}}(g_{1}(\mathbf{x})|D_{1}%
u|^{2}+\varepsilon^{2}g_{2}(\mathbf{x})|D_{2}u|^{2})d\mathbf{x},\ \ u\in
H^{1}({\mathbb{R}}^{2}). \label{2.1}%
\end{equation}
We have the obvious identity%
\begin{equation}
A_{\varepsilon}=\varepsilon^{-2}T_{\varepsilon}^{\ast}A(\varepsilon
)T_{\varepsilon}. \label{2.2}%
\end{equation}
Clearly, for the operator $\left[  Q^{\varepsilon}\right]  $ of multiplication
by the function $Q^{\varepsilon}(\mathbf{x})$, we have%
\begin{equation}
\left[  Q^{\varepsilon}\right]  =T_{\varepsilon}^{\ast}\left[  Q\right]
T_{\varepsilon}. \label{2.3}%
\end{equation}
From (\ref{2.2}) and (\ref{2.3}) it follows that%
\begin{equation}
(A_{\varepsilon}+Q^{\varepsilon})^{-1}=\varepsilon^{2}T_{\varepsilon}^{\ast
}(A(\varepsilon)+\varepsilon^{2}Q)^{-1}T_{\varepsilon}. \label{2.4}%
\end{equation}
A similar representation is true for $(A^{0}+Q^{0})^{-1}$:%
\begin{equation}
(A^{0}+Q^{0})^{-1}=\varepsilon^{2}T_{\varepsilon}^{\ast}(A^{0}(\varepsilon
)+\varepsilon^{2}Q^{0})^{-1}T_{\varepsilon}, \label{2.5}%
\end{equation}
where $A^{0}(\varepsilon)$ corresponds to the quadratic form%
\begin{equation}
a^{0}(\varepsilon)[u,u]=\int_{{\mathbb{R}}^{2}}(g_{1}^{0}(x_{2})|D_{1}%
u|^{2}+\varepsilon^{2}g_{2}^{0}(x_{2})|D_{2}u|^{2})\,d\mathbf{x},\ \ u\in
H^{1}({\mathbb{R}}^{2}). \label{2.6}%
\end{equation}
We can also define $A^{0}(\varepsilon)$ by the differential expression
$g_{1}^{0}(x_{2})D_{1}^{2}+\varepsilon^{2}D_{2}g_{2}^{0}(x_{2})D_{2}$ on
domain $H^{2}({\mathbb{R}}^{2}).$

Using (\ref{2.4}), (\ref{2.5}) and the fact that $T_{\varepsilon}$ is unitary
operator in $L_{2}({\mathbb{R}}^{2})$, we reduce Theorem \ref{theor1.1} to the
following theorem.

\begin{theorem}
\label{theor2.1}Suppose that conditions \textrm{(\ref{1.1})-(\ref{1.3})} are
satisfied. Let $A(\varepsilon)$ be the operator in $L_{2}({\mathbb{R}}^{2})$
that corresponds to the quadratic form \textrm{(\ref{2.1})}. Let
$A^{0}(\varepsilon)$ be the operator corresponding to the form
\textrm{(\ref{2.6})} with coefficients defined by \textrm{(\ref{1.11}),
(\ref{1.12})}. Let $Q(\mathbf{x})$ be a function in ${\mathbb{R}}^{2}$
satisfying \textrm{(\ref{1.6a})--(\ref{1.8})}, and let $Q^{0}(x_{2})$ be
defined by \textrm{(\ref{1.15})}. Then we have%
\begin{equation}
\|(A(\varepsilon)+\varepsilon^{2}Q)^{-1}-(A^{0}(\varepsilon)+\varepsilon
^{2}Q^{0})^{-1}\|_{L_{2}({\mathbb{R}}^{2})\rightarrow L_{2}({\mathbb{R}}^{2}%
)}\leq C\varepsilon^{-1},\ \ 0<\varepsilon\leq1. \label{2.7}%
\end{equation}
The constant $C$ depends only on $c_{j},\,j=0,...,5.$
\end{theorem}

\subsection{Direct integral decomposition for $A(\varepsilon)\label{sect2.2}$}

First, we define the Gelfand transformation $\mathcal{{U}}$. We denote%
\begin{equation}
\Omega=(0,1)\times{\mathbb{R}}, \label{2.8}%
\end{equation}
and consider the Hilbert space%
\begin{equation}
\mathcal{K}=\int_{\left[  -\pi,\pi\right)  }\oplus L_{2}(\Omega)dk=L_{2}%
([-\pi,\pi);L_{2}(\Omega))=L_{2}([-\pi,\pi)\times\Omega), \label{2.9}%
\end{equation}
which is the direct integral with constant fibers. Initially, the operator
$\mathcal{{U}}:L_{2}\mathcal{({\mathbb{R}}}^{2}\mathcal{)}\rightarrow
\mathcal{{K}}$ is defined on the functions $f\in C_{0}^{\infty}({\mathbb{R}%
}^{2})$ by the formula%
\begin{equation}
(\mathcal{U}f)(\mathbf{x},k)=(2\pi)^{-1/2}\sum_{m\in{\mathbb{Z}}}%
e^{-ik(x_{1}+m)}f(x_{1}+m,x_{2}),\ \ \mathbf{x}=(x_{1},x_{2})\in\Omega
,\;k\in\lbrack-\pi,\pi). \label{2.10}%
\end{equation}
Next, $\mathcal{U}$ extends by continuity to a \textit{unitary mapping of
$L_{2}({\mathbb{R}}^{2})$ onto} $\mathcal{K}$.

In $L_{2}(\Omega)$, we consider the operator $A(k,\varepsilon)$ generated by
the quadratic form%
\begin{equation}
a(k,\varepsilon)[u,u]=\int_{\Omega}(g_{1}(\mathbf{x})|(D_{1}+k)u|^{2}%
+\varepsilon^{2}g_{2}(\mathbf{x})|D_{2}u|^{2})\,d\mathbf{x},\quad
u\in\widetilde{H}^{1}(\Omega),\ k\in\lbrack-\pi,\pi),\ \varepsilon>0.
\label{2.11}%
\end{equation}
Here by $\widetilde{H}^{s}(\Omega)$ we denote the subspace of $H^{s}\left(
\Omega\right)  $ formed by the functions whose 1-periodic extension (in
$x_{1})$ to $\mathbb{R}^{2}$ belongs to $H_{\mathrm{loc}}^{s}(\mathbb{R}%
^{2}).$

\textit{The operator $A(\varepsilon)$ is decomposed in the direct integral of
operators} $A(k,\varepsilon)$ with the help of the Gelfand transformation
$\mathcal{U}$:%
\[
\mathcal{U}A(\varepsilon)\mathcal{U}^{-1}=\int_{[-\pi,\pi)}\oplus
A(k,\varepsilon)dk.
\]
Then, for the generalized resolvent $(A(\varepsilon)+\varepsilon^{2}Q)^{-1}$
we have%
\begin{equation}
\mathcal{U}(A(\varepsilon)+\varepsilon^{2}Q)^{-1}\mathcal{U}^{-1}=\int
_{[-\pi,\pi)}\oplus(A(k,\varepsilon)+\varepsilon^{2}Q)^{-1}dk. \label{2.12}%
\end{equation}

Let $A^{0}(k,\varepsilon)$ be the operator in $L_{2}(\Omega)$ corresponding to
the quadratic form%
\begin{equation}
a^{0}(k,\varepsilon)[u,u]=\int_{\Omega}(g_{1}^{0}(x_{2})|(D_{1}+k)u|^{2}%
+\varepsilon^{2}g_{2}^{0}(x_{2})|D_{2}u|^{2})\,d\mathbf{x},\ \ u\in
\widetilde{H}^{1}(\Omega),\,\,k\in\lbrack-\pi,\pi),\,\varepsilon>0.
\label{2.13}%
\end{equation}
We can also define $A^{0}(k,\varepsilon)$ by the differential expression
$g_{1}^{0}(x_{2})(D_{1}+k)^{2}+\varepsilon^{2}D_{2}g_{2}^{0}(x_{2})D_{2}$ on
the domain $\widetilde{H}^{2}(\Omega)$. For $A^{0}(\varepsilon)$ we also have
the direct integral decomposition:%
\begin{equation}
\mathcal{U}A^{0}(\varepsilon)\mathcal{U}^{-1}=\int_{[-\pi,\pi)}\oplus
A^{0}(k,\varepsilon)dk. \label{2.14}%
\end{equation}
Then, for the generalized resolvent $(A^{0}(\varepsilon)+\varepsilon^{2}%
Q^{0})^{-1}$ we have%
\begin{equation}
\mathcal{U}(A^{0}(\varepsilon)+\varepsilon^{2}Q^{0})^{-1}\mathcal{U}^{-1}%
=\int_{[-\pi,\pi)}\oplus(A^{0}(k,\varepsilon)+\varepsilon^{2}Q^{0})^{-1}dk.
\label{2.15}%
\end{equation}

Using (\ref{2.12}), (\ref{2.15}) and the fact that $\mathcal{U}$ is unitary
transformation, we see that Theorem \ref{theor2.1} (and then also Theorem
\ref{theor1.1}) is a direct consequence of the following statement.

\begin{theorem}
\label{theor2.2}Suppose that conditions \textrm{(\ref{1.1})--(\ref{1.3})} are
satisfied. Let $A(k,\varepsilon)$ be the operator in $L_{2}(\Omega)$ generated
by the quadratic form \textrm{(\ref{2.11})}. Let $A^{0}(k,\varepsilon)$ be the
operator in $L_{2}(\Omega)$ corresponding to the form \textrm{(\ref{2.13})},
where the coefficients are defined by \textrm{(\ref{1.11}), (\ref{1.12})}. Let
$Q(\mathbf{x})$ be a function satisfying \textrm{(\ref{1.6a})--(\ref{1.8})},
and let $Q^{0}(x_{2})$ be defined by \textrm{(\ref{1.15})}. Then we have%
\begin{equation}
\left\Vert (A(k,\varepsilon)+\varepsilon^{2}Q)^{-1}-(A^{0}(k,\varepsilon
)+\varepsilon^{2}Q^{0})^{-1}\right\Vert _{L_{2}(\Omega)\rightarrow
L_{2}(\Omega)}\leq C\varepsilon^{-1},\quad\ 0<\varepsilon\leq1,\ k\in
\lbrack-\pi,\pi). \label{2.16}%
\end{equation}
The constant $C$ depends only on $c_{j},j=0,...,5.$
\end{theorem}

\section{The operator pencil $A_{1}(k)$\label{sect3}}

\subsection{Definition of $A_{1}(k)$. Direct integral decomposition}

If we formally put $\varepsilon=0$, then the operator $A(k,\varepsilon)$
defined in Subsection \ref{sect2.2} turns into the operator $A_{1}(k)$
formally given by the expression $(D_{1}+k)g_{1}(\mathbf{x})(D_{1}+k)$ with
periodic boundary conditions on $\partial\Omega$. Precisely, $A_{1}(k)$ is the
selfadjoint operator in $L_{2}(\Omega)$ generated by the quadratic form%
\begin{equation}
a_{1}(k)[u,u]=\int_{\Omega}g_{1}(\mathbf{x})|(D_{1}+k)u|^{2}\,d\mathbf{x},
\label{3.1}%
\end{equation}%
\begin{equation}
\operatorname*{Dom}a_{1}(k)=L_{2}({\mathbb{R}};\widetilde{H}^{1}%
(0,1))=\int_{\mathbb{R}}\oplus\widetilde{H}^{1}(0,1)dx_{2}. \label{3.2}%
\end{equation}
Here $\widetilde{H}^{1}(0,1)$ is the subspace of $H^{1}(0,1)$ formed by the
functions whose 1-periodic extension to $\mathbb{R}$ belongs to
$H_{\mathrm{loc}}^{1}(\mathbb{R}).$

In this section we study the operator family $A_{1}(k)$. For this, it suffices
to impose only condition (\ref{1.1}) with $j=1$. Now it is convenient to
interpret the space $L_{2}(\Omega)$ as a direct integral in $x_{2}$ with
constant fibers $L_{2}(0,1)$:%
\begin{equation}
L_{2}(\Omega)= \int_{\mathbb{R}}\oplus L_{2}(0,1)dx_{2}. \label{3.3}%
\end{equation}
The operator $A_{1}(k)$ acts in the space (\ref{3.3}) \textquotedblleft
layerwise\textquotedblright, i.~e., it is the operator of multiplication by
the operator-valued function $A_{1}(k;x_{2})$:%
\begin{equation}
A_{1}(k)=\int_{\mathbb{R}}\oplus A_{1}(k;x_{2})dx_{2}, \label{3.4}%
\end{equation}
where $A_{1}(k,x_{2})$ is the selfadjoint operator in $L_{2}(0,1)$ generated
by the form%
\begin{equation}
a_{1}(k;x_{2})[u,u]=\int_{0}^{1}g_{1}(x_{1},x_{2})|(D_{1}+k)u|^{2}%
dx_{1},\,\ \ u\in\widetilde{H}^{1}(0,1),\ x_{2}\in{\mathbb{R}},\ k\in
\lbrack-\pi,\pi). \label{3.5}%
\end{equation}

\subsection{The operator family $A_{1}(k;x_{2})$}

The operator family $A_{1}(k;x_{2})$ can be studied by the method suggested in
\cite{BSu1}. This family admits a factorization%
\begin{equation}
A_{1}(k;x_{2})=\left(  X(k;x_{2})\right)  ^{\ast}X(k;x_{2}), \label{3.6}%
\end{equation}%
\begin{equation}
X(k;x_{2})=X_{0}(x_{2})+kX_{1}(x_{2}), \label{3.7}%
\end{equation}
where $X_{0}(x_{2})$ is the (closed) operator in $L_{2}(0,1)$ given by the
differential expression $(g_{1}(x_{1},x_{2}))^{\frac{1}{2}}D_{1}$ on the
domain $\widetilde{H}^{1}(0,1)$, and $X_{1}(x_{2})$ is the bounded operator in
$L_{2}(0,1)$ of multiplication by the function $(g_{1}(x_{1},x_{2}))^{\frac
{1}{2}}$. Then the operator (\ref{3.7}) is closed on the domain $\widetilde
{H}^{1}(0,1)$. Obviously, the kernel
\begin{equation}
\widetilde{\mathfrak{N}}=\operatorname*{Ker}A_{1}(0;x_{2})=\operatorname*{Ker}%
X_{0}(x_{2})=\{{u\in L_{2}(0,1):u=}\text{const}{\}} \label{3.8}%
\end{equation}
is one-dimensional. Let $d_{0}(x_{2})$ be the distance from the point
$\lambda_{0}=0$ to the rest of the spectrum of $A_{1}(0;x_{2})$. Let us
estimate $d_{0}(x_{2})$ from below.

By (\ref{1.1}) (with $j=1$), the form (\ref{3.5}) with $k=0$ satisfies the
estimate%
\begin{equation}
a_{1}(0;x_{2})[u,u]\geq c_{0}\int_{0}^{1}|D_{1}u(x_{1})|^{2}dx_{1}%
,\;\ u\in{\widetilde{H}}^{1}(0,1). \label{3.9}%
\end{equation}
Using the Fourier series expansion for $u\in{\widetilde{H}}^{1}(0,1)$
\[
u(x_{1})=\sum_{m\in\mathbb{Z}}\widehat{u}_{m}e^{2\pi imx_{1}},
\]
and (\ref{3.9}), we obtain%
\begin{equation}
a_{1}(0;x_{2})[u,u]\geq c_{0}\sum_{m\in\mathbb{Z}}4\pi^{2}m^{2}|\widehat
{u}_{m}|^{2},\ \;u\in{\widetilde{H}}^{1}(0,1). \label{3.10}%
\end{equation}
If $u$ is orthogonal to $\widetilde{\mathfrak{N}}$ in $L_{2}(0,1)$, then
$\widehat{u}_{0}=0$, and from (\ref{3.10}) we deduce that%
\begin{equation}
a_{1}(0;x_{2})[u,u]\geq4\pi^{2}c_{0}\|u\|_{L_{2}(0,1)}^{2},\quad
u\in{\widetilde{H}}^{1}(0,1),\ u\perp\widetilde{\mathfrak{N}}. \label{3.11}%
\end{equation}
It follows that%
\begin{equation}
d_{0}(x_{2})\geq4\pi^{2}c_{0}=:d_{0},\ \ x_{2}\in\mathbb{R}. \label{3.12}%
\end{equation}
In accordance with \cite[Ch.1, \S 1]{BSu1}, we choose the positive number
$\delta$ such that $8 \delta< d_{0}.$ We can put
\begin{equation}
\delta=\frac{\pi^{2}c_{0}}{4}. \label{3.13}%
\end{equation}
Next, by the upper estimate (\ref{1.1}), we have%
\begin{equation}
\|X_{1}(x_{2})\|_{L_{2}(0,1)\rightarrow L_{2}(0,1)}\leq\|g_{1}(\cdot
,x_{2})\|_{L_{\infty}(0,1)}^{\frac{1}{2}}\leq c_{1}^{\frac{1}{2}},\ \;x_{2}%
\in\mathbb{R}. \label{3.14}%
\end{equation}
As in \cite[Ch. 1, \S 1]{BSu1}, we choose a number $t_{0}$ such that
$t_{0}\leq\delta^{\frac{1}{2}}\|X_{1}(x_{2})\|^{-1}$ for all $x_{2}%
\in\mathbb{R}.$ We put
\begin{equation}
t_{0}={\frac{\pi}{2}} c_{0}^{{\frac{1}{2}}} c_{1}^{-{\frac{1}{2}}}.
\label{3.15}%
\end{equation}

Let $F(k;x_{2};\sigma)$ be the spectral projection of the operator
$A_{1}(k;x_{2})$ corresponding to the interval $[0,\sigma]$. From Proposition
1.2 of \cite[Ch. 1]{BSu1} it follows that%
\begin{equation}
\dim F(k;x_{2};\delta)L_{2}(0,1)=\dim F(k;x_{2};3\delta)L_{2}%
(0,1)=1,\ \ \,|k|=t\leq t_{0}. \label{3.16}%
\end{equation}
It means that for $|k|\leq t_{0}$ and all $x_{2}\in\mathbb{R}$ the operator
$A_{1}(k;x_{2})$ has only one eigenvalue lying in the interval $[0,\delta]$,
while the interval $(\delta,3\delta)$ is free of its spectrum. We put%
\begin{equation}
F(k;x_{2})=F(k;x_{2};\delta). \label{3.17}%
\end{equation}
Let $\widetilde{P}$ be the orthogonal projection of $L_{2}(0,1)$ onto the
subspace $\widetilde{\mathfrak{N}}$. By Theorem 4.1 of \cite[Ch.1]{BSu1}, we
have the following approximation for the operator (\ref{3.17}):%
\begin{equation}
F(k;x_{2})=\widetilde{P}+t\Phi(k;x_{2}),\quad t=|k|\leq t_{0}, \ x_{2}%
\in\mathbb{R}, \label{3.18}%
\end{equation}%
\begin{equation}
\|\Phi(k;x_{2})\|_{L_{2}(0,1)\rightarrow L_{2}(0,1)}\leq C_{1},\quad t\leq
t_{0},\;x_{2}\in\mathbb{R}, \label{3.19}%
\end{equation}
where
\begin{equation}
C_{1}=\beta_{1}\delta^{-\frac{1}{2}}\sup_{x_{2}\in\mathbb{R}}\|X_{1}%
(x_{2})\|\leq2\beta_{1}c_{1}^{\frac{1}{2}}c_{0}^{-\frac{1}{2}}\pi^{-1},
\label{3.20}%
\end{equation}
and $\beta_{1}$ is an absolute constant.

In accordance with \cite{BSu1}, we introduce the spectral germ $S(x_{2})$ of
the operator family $A_{1}(k;x_{2})$. The germ is a selfadjoint operator
acting in the space $\widetilde{\mathfrak{N}}$. In our case (i.~e., for the
operator of the form $-\frac{d}{dx}g(x)\frac{d}{dx}$), the germ is calculated
explicitly (see \cite[Ch.5, \S 1]{BSu1}). The operator $S(x_{2})$ is
multiplication by the number $g_{1}^{0}(x_{2})$, defined by (\ref{1.11}).

By Theorem 4.3 of \cite[Ch.1]{BSu1}, we have the following approximation for
the operator $A_{1}(k;x_{2})F(k;x_{2})$:%
\begin{equation}
A_{1}(k;x_{2})F(k;x_{2})=t^{2}S(x_{2})\widetilde{P}+t^{3}\Psi(k;x_{2}%
),\;\ t=|k|\leq t_{0},\,x_{2}\in\mathbb{R}, \label{3.21}%
\end{equation}%
\begin{equation}
\|\Psi(k;x_{2})\|_{L_{2}(0,1)\rightarrow L_{2}(0,1)}\leq C_{2},\ \ t\leq
t_{0},\,x_{2}\in\mathbb{R}, \label{3.22}%
\end{equation}
where%
\begin{equation}
C_{2}=\beta_{2}\delta^{-\frac{1}{2}}\sup_{x_{2}\in\mathbb{R}}\|X_{1}%
(x_{2})\|^{3}\leq2\beta_{2}\pi^{-1}c_{1}^{\frac{3}{2}}c_{0}^{-\frac{1}{2}},
\label{3.23}%
\end{equation}
and $\beta_{2}$ is an absolute constant.

\subsection{The spectral germ of the operator family $A_{1}(k)$}

We return to the family (\ref{3.4}) acting in $L_{2}(\Omega)$. From
(\ref{3.4}) and (\ref{3.16}) it follows that the interval $(\delta,3\delta)$
is free of the spectrum of $A_{1}(k)$, if $|k|\leq t_{0}.$ By (\ref{3.4}) and
(\ref{3.6}), we have $A_{1}(k)=X(k)^{\ast}X(k)$, where $X(k)=\int_{\mathbb{R}%
}\oplus X(k;x_{2})dx_{2}.$ Then, by (\ref{3.7}), $X(k)=X_{0}+kX_{1}%
,\,X_{0}=\int_{\mathbb{R}}\oplus X_{0}(x_{2})dx_{2},\,X_{1}=\int_{\mathbb{R}%
}\oplus X_{1}(x_{2})dx_{2}.$

Clearly, $X_{0}$ is the operator in $L_{2}(\Omega)$ given by $(g_{1}%
(x))^{\frac{1}{2}}D_{1}$ on the domain (\ref{3.2}), and $X_{1}$ is the bounded
operator in $L_{2}(\Omega)$ of multiplication by the function $(g_{1}%
(x))^{\frac{1}{2}}$.

For the operator $A_{1}(0)$ with $k=0$ we have%
\begin{equation}
\mathfrak{N}=\operatorname*{Ker}A_{1}(0)=\operatorname*{Ker}X_{0}=\{{u\in
L_{2}(\Omega):u=u(x_{2})\}}. \label{3.24}%
\end{equation}
Thus, the kernel $\mathfrak{N}$ consists of functions in $L_{2}(\Omega)$
depending only on the second variable $x_{2}$. It can be identified with
$L_{2}(\mathbb{R})$. Obviously, we can write%
\begin{equation}
\mathfrak{N}=\int_{\mathbb{R}}\oplus\widetilde{\mathfrak{N}}\,dx_{2}.
\label{3.24b}%
\end{equation}
Also, the orthogonal projection $P$ of $L_{2}(\Omega)$ onto $\mathfrak{N}$ can
be represented as the direct integral%
\begin{equation}
P=\int_{\mathbb{R}}\oplus\widetilde{P}\,dx_{2}. \label{3.25}%
\end{equation}
Now, we consider the operator $S$ in $\mathfrak{N}$ which is the direct
integral of $S(x_{2})$:
\begin{equation}
S\overset{\mathrm{def}}{=}\int_{\mathbb{R}}\oplus S(x_{2})\,dx_{2}.
\label{3.26}%
\end{equation}
Then $S$ is the operator in $\mathfrak{N}\equiv L_{2}(\mathbb{R})$ of
multiplication by the function $g_{1}^{0}(x_{2}).$ \textit{The operator $S$ is
called the spectral germ of $A_{1}(k).$}

Let $F(k)$ denote the spectral projection of the operator $A_{1}(k)$
corresponding to the interval $[0,\delta]$. Since $A_{1}(k)$ is represented as
the direct integral (\ref{3.4}), we have%
\begin{equation}
F(k)=\int_{\mathbb{R}}\oplus F(k;x_{2})\,dx_{2}. \label{3.27}%
\end{equation}
From (\ref{3.18}), (\ref{3.25}) and (\ref{3.27}) it follows that%
\begin{equation}
F(k)=P+t\Phi(k),\quad t=|k|\leq t_{0}, \label{3.28}%
\end{equation}%
\begin{equation}
\Phi(k)=\int_{\mathbb{R}}\oplus\Phi(k;x_{2})\,dx_{2}. \label{3.29}%
\end{equation}
Similarly, relations (\ref{3.4}), (\ref{3.21}), (\ref{3.25})-(\ref{3.27})
imply that%
\begin{equation}
A_{1}(k)F(k)=t^{2}SP+t^{3}\Psi(k),\quad t=|k|\leq t_{0}, \label{3.30}%
\end{equation}%
\begin{equation}
\Psi(k)=\int_{\mathbb{R}}\oplus\Psi(k;x_{2})\,dx_{2}. \label{3.31}%
\end{equation}
By (\ref{3.19}) and (\ref{3.22}), the operators (\ref{3.29}) and (\ref{3.31})
satisfy the estimates%
\begin{equation}
\|\Phi(k)\|_{L_{2}(\Omega)\rightarrow L_{2}(\Omega)}\leq C_{1}, \quad
t=|k|\leq t_{0}, \label{3.32}%
\end{equation}%
\begin{equation}
\|\Psi(k)\|_{L_{2}(\Omega)\rightarrow L_{2}(\Omega)}\leq C_{2}, \quad
t=|k|\leq t_{0}. \label{3.33}%
\end{equation}

\section{Approximation for the generalized resolvent of $A\left(
k,\varepsilon\right)  $\label{sect4}}

\subsection{Approximation in terms of the operator $S(k,\varepsilon)$}

Now we return to the operator family $A(k,\varepsilon)$ introduced in
Subsection \ref{sect2.2}. Our goal is to approximate the generalized resolvent
$(A(k,\varepsilon)+\varepsilon^{2}Q)^{-1}.$

In the space $\mathfrak{N}\equiv L_{2}(\mathbb{R})$ (see (\ref{3.24})), we
consider the quadratic form%
\begin{equation}
s(k,\varepsilon)[\omega,\omega]=\int_{\mathbb{R}}(g_{1}^{0}(x_{2})k^{2}%
|\omega(x_{2})|^{2}+\varepsilon^{2}g_{2}^{0}(x_{2})|D_{2}\omega(x_{2}%
)|^{2})dx_{2},\ \ k\in\lbrack-\pi,\pi),\ \varepsilon>0, \label{4.1}%
\end{equation}%
\begin{equation}
\operatorname*{Dom}s(k,\varepsilon)=H^{1}(\mathbb{R}). \label{4.2}%
\end{equation}
The selfadjoint operator in $L_{2}(\mathbb{R})$ corresponding to this form is
denoted by $S(k,\varepsilon).$ Clearly, under conditions (\ref{1.13}),
(\ref{1.14a}), the operator $S(k,\varepsilon)$ is given by the differential
expression%
\begin{equation}
S(k,\varepsilon)=g_{1}^{0}(x_{2})k^{2}+\varepsilon^{2}D_{2}g_{2}^{0}%
(x_{2})D_{2}, \label{4.3}%
\end{equation}
on the domain%
\begin{equation}
\operatorname*{Dom}S(k,\varepsilon)=H^{2}(\mathbb{R}). \label{4.4}%
\end{equation}
Note that the first summand in (\ref{4.3}) corresponds to $t^{2}S.$

The key role in the proof of Theorem \ref{theor2.2} is played by the following
theorem, which gives an approximation for the generalized resolvent of
$A(k,\varepsilon)$ in terms of the operator $S(k,\varepsilon)$.

\begin{theorem}
\label{theor4.1}Suppose that conditions \textrm{(\ref{1.1})--(\ref{1.3})} are
satisfied. Let $A(k,\varepsilon)$ be the operator in $L_{2}(\Omega)$ generated
by the quadratic form \textrm{(\ref{2.11})}. Let $S(k,\varepsilon)$ be the
operator in $\mathfrak{N}$ defined by \textrm{(\ref{4.3}), (\ref{4.4})}, where
the coefficients are defined by \textrm{(\ref{1.11}), (\ref{1.12})}. Let
$Q(\mathbf{x})$ be a function satisfying \textrm{(\ref{1.6a})--(\ref{1.8})}
and let $Q^{0}(x_{2})$ be defined by \textrm{(\ref{1.15})}. Let $P$ be the
orthogonal projection of $L_{2}(\Omega)$ onto the subspace \textrm{(\ref{3.24}%
)}. Then we have%
\begin{equation}
\|(A(k,\varepsilon)+\varepsilon^{2}Q)^{-1}-(S(k,\varepsilon)+\varepsilon
^{2}Q^{0})^{-1}P\|_{L_{2}(\Omega)\rightarrow L_{2}(\Omega)}\leq\widetilde
{C}\varepsilon^{-1},\ \ |k|\leq t_{0},\ 0<\varepsilon\leq1. \label{4.5}%
\end{equation}
The number $t_{0}$ is defined by \textrm{(\ref{3.15})}. The constant
$\widetilde{C}$ depends only on $c_{j},\,j=0,...,5.$
\end{theorem}

\subsection{Deduction of Theorem \ref{theor2.2} from Theorem \ref{theor4.1}%
\label{sect4.2}$_{{}}$}

Let $A^{0}(k,\varepsilon)$ be the operator in $L_{2}(\Omega)$ corresponding to
the form (\ref{2.13}). For $A^{0}(k,\varepsilon)$, we define the operator
$S^{0}(k,\varepsilon)$ by the same rule as for $A(k,\varepsilon).$ Clearly,
$S^{0}(k,\varepsilon)=S(k,\varepsilon).$ Applying Theorem \ref{theor4.1} to
the operator $(A^{0}(k,\varepsilon)+ \varepsilon^{2} Q^{0})^{-1}$, we obtain%
\begin{equation}
\|(A^{0}(k,\varepsilon)+\varepsilon^{2}Q^{0})^{-1}-(S(k,\varepsilon
)+\varepsilon^{2}Q^{0})^{-1}P\|_{L_{2}(\Omega)\rightarrow L_{2}(\Omega)}%
\leq\widetilde{C}^{\prime} \varepsilon^{-1},\ \ |k|\leq t_{0},\ 0<\varepsilon
\leq1. \label{4.6}%
\end{equation}
Since the coefficients $g_{j}^{0}(x_{2}),\,j=1,2$, satisfy estimates
(\ref{1.13}) with the same constants $c_{0}$ and $c_{1}$ as in (\ref{1.1}),
the number $t_{0}$ in (\ref{4.5}) and (\ref{4.6}) is one and the same (see
(\ref{3.15})). Also, using the statement of Theorem \ref{theor4.1} concerning
the constant $\widetilde{C}$ and (\ref{1.14}), (\ref{1.14a}), we see that
$\widetilde{C}^{\prime}$ is controlled in terms of the constants
$c_{j},j=0,...,5.$

Combining (\ref{4.5}) and (\ref{4.6}), we have%
\begin{equation}
\|(A(k,\varepsilon)+\varepsilon^{2}Q)^{-1}-(A^{0}(k,\varepsilon)+\varepsilon
^{2}Q^{0})^{-1}\|_{L_{2}(\Omega)\rightarrow L_{2}(\Omega)}\leq(\widetilde
{C}+\widetilde{C}^{\prime})\varepsilon^{-1},\quad|k|\leq t_{0},\ 0<\varepsilon
\leq1. \label{4.7}%
\end{equation}
So, we have obtained the required estimate (\ref{2.16}) for $|k|\leq t_{0}$.

If $k\in\lbrack-\pi,\pi)$ and $|k|>t_{0}$, the estimates are trivial: both
terms in (\ref{2.16}) are estimated by constants. Indeed, by (\ref{1.1}) and
(\ref{2.11}),%
\[
a(k,\varepsilon)[u,u]\geq c_{0}\int_{\Omega}(|(D_{1}+k)u|^{2}+\varepsilon
^{2}|D_{2}u|^{2})\,d\mathbf{x},\quad u\in\widetilde{H}^{1}(\Omega
),\ k\in\lbrack-\pi,\pi),\ \varepsilon>0.
\]
Using the Fourier series expansion (in variable $x_{1}$) for $u\in
\widetilde{H}^{1}(\Omega)$:%
\[
u(x_{1},x_{2})= \sum_{m\,\in\,\mathbb{Z}}\widehat{u}_{m}(x_{2})e^{2\pi
imx_{1}},
\]
we obtain%
\begin{equation}
a(k,\varepsilon)[u,u]\geq c_{0}\sum_{m\in\mathbb{Z}}(2\pi m+k)^{2}%
\int_{\mathbb{R}}|\widehat{u}_{m}(x_{2})|^{2}dx_{2},\quad u\in\widetilde
{H}^{1}(\Omega),\ k\in\lbrack-\pi,\pi),\ \varepsilon>0. \label{4.8}%
\end{equation}
If $k\in\lbrack-\pi,-t_{0})\cup(t_{0},\pi)$, then $|2\pi m+k|\geq t_{0}%
,\,m\in\mathbb{Z}.$ Then (\ref{4.8}) implies that%
\begin{equation}
a(k,\varepsilon)[u,u]\geq c_{0}t_{0}^{2}\|u\|_{L_{2}(\Omega)}^{2},\quad
u\in\widetilde{H}^{1}(\Omega),\ k\in\lbrack-\pi,-t_{0})\cup(t_{0}%
,\pi),\ \varepsilon>0. \label{4.9}%
\end{equation}
Consequently,%
\begin{equation}
A(k,\varepsilon)\geq c_{0}t_{0}^{2}I,\quad k\in\lbrack-\pi,-t_{0})\cup
(t_{0},\pi),\ \varepsilon>0, \label{4.10}%
\end{equation}
and%
\begin{equation}
\|(A(k,\varepsilon)+\varepsilon^{2}Q)^{-1}\|_{L_{2}(\Omega)\rightarrow
L_{2}(\Omega)}\leq c_{0}^{-1}t_{0}^{-2},\quad k\in\lbrack-\pi,-t_{0}%
)\cup(t_{0},\pi),\ \varepsilon>0. \label{4.11}%
\end{equation}
The same estimate is true for $(A^{0}(k,\varepsilon)+\varepsilon^{2}%
Q^{0})^{-1}.$ Hence,%
\begin{equation}
\|(A(k,\varepsilon)+\varepsilon^{2}Q)^{-1}-(A^{0}(k,\varepsilon)+\varepsilon
^{2}Q^{0})^{-1}\|_{L_{2}(\Omega)\rightarrow L_{2}(\Omega)}\leq2c_{0}^{-1}%
t_{0}^{-2},\quad k\in\lbrack-\pi,-t_{0})\cup(t_{0},\pi),\ \varepsilon>0.
\label{4.12}%
\end{equation}
From (\ref{4.7}) and (\ref{4.12}) it follows that (\ref{2.16}) is valid with
the constant $C=\max\{{\widetilde{C}+\widetilde{C}^{\prime},\,2c_{0}^{-1}%
t_{0}^{-2}\}}.$

Theorem \ref{theor4.1} will be proved in Sections \ref{sect5}, \ref{sect6},
\ref{sect7}. The proof needs some technical preparations.

\section{Preliminary estimates\label{sect5}}

\subsection{Estimates for generalized resolvents of the operators
$A(k,\varepsilon)$, $A_{1}(k)$ and $S(k,\varepsilon).$}

\begin{lemma}
\label{lem5.1}We have%
\begin{equation}
\|(A(k,\varepsilon)+\varepsilon^{2}Q)^{-1}\|_{L_{2}(\Omega)\rightarrow
L_{2}(\Omega)}\leq(c_{0}k^{2}+c_{3}\varepsilon^{2})^{-1},\quad k\in\lbrack
-\pi,\pi),\ \varepsilon>0, \label{5.1}%
\end{equation}%
\begin{equation}
\|(A_{1}(k)+\varepsilon^{2}Q)^{-1}\|_{L_{2}(\Omega)\rightarrow L_{2}(\Omega
)}\leq(c_{0}k^{2}+c_{3}\varepsilon^{2})^{-1},\quad k\in\lbrack-\pi
,\pi),\ \varepsilon>0, \label{5.2}%
\end{equation}%
\begin{equation}
\|(S(k,\varepsilon)+\varepsilon^{2}Q^{0})^{-1}P\|_{L_{2}(\Omega)\rightarrow
L_{2}(\Omega)}\leq(c_{0}k^{2}+c_{3}\varepsilon^{2})^{-1},\quad k\in\lbrack
-\pi,\pi),\ \varepsilon>0. \label{5.3}%
\end{equation}

\end{lemma}

\textbf{Proof.} From (\ref{4.8}) it follows that%
\begin{equation}
a(k,\varepsilon)[u,u]\geq c_{0}k^{2}\|u\|_{L_{2}(\Omega)}^{2},\quad
u\in\widetilde{H}^{1}(\Omega),\ k\in\lbrack-\pi,\pi),\ \varepsilon>0.
\label{5.4}%
\end{equation}
Combining this with the lower estimate (\ref{1.7}), we see that%
\begin{equation}
a(k,\varepsilon)[u,u]+\varepsilon^{2}\int_{\Omega}Q(\mathbf{x})|u|^{2}%
\,d\mathbf{x}\geq(c_{0}k^{2}+c_{3}\varepsilon^{2})\|u\|_{L_{2}(\Omega)}%
^{2},\ \ u\in\widetilde{H}^{1}(\Omega),\,k\in\lbrack-\pi,\pi),\,\varepsilon>0,
\label{5.5}%
\end{equation}
which implies (\ref{5.1}).

From (\ref{1.1}) and (\ref{3.1}) it follows that%
\begin{equation}
a_{1}(k)[u,u]\geq c_{0}\int_{\Omega}|(D_{1}+k)u|^{2}\,d\mathbf{x},\quad
u\in\widetilde{H}^{1}(\Omega),\ k\in\lbrack-\pi,\pi),\ \varepsilon>0.
\label{5.6}%
\end{equation}
Hence, similarly to (\ref{4.8}),%
\begin{equation}
a_{1}(k)[u,u]\geq c_{0}k^{2}\|u\|_{L_{2}(\Omega)}^{2},\quad u\in\widetilde
{H}^{1}(\Omega),\ k\in\lbrack-\pi,\pi),\ \varepsilon>0. \label{5.7}%
\end{equation}
Relations (\ref{5.7}) and (\ref{1.7}) imply that%
\begin{equation}
a_{1}(k)[u,u]+\varepsilon^{2}\int_{\Omega}Q(\mathbf{x})|u|^{2}\,d\mathbf{x}%
\geq(c_{0}k^{2}+c_{3}\varepsilon^{2})\|u\|_{L_{2}(\Omega)}^{2},\quad
u\in\widetilde{H}^{1}(\Omega),\ k\in\lbrack-\pi,\pi),\ \varepsilon>0.
\label{5.8}%
\end{equation}
This yields (\ref{5.2}).

Finally, by (\ref{4.1}) and (\ref{1.13}),%
\begin{align}
s(k,\varepsilon)[\omega,\omega]  &  \geq c_{0}\int_{\mathbb{R}}(k^{2}%
|\omega(x_{2})|^{2}+\varepsilon^{2}|D_{2}\omega(x_{2})|^{2})dx_{2}%
\geq\label{5.9}\\
&  \geq c_{0}k^{2}\|\omega\|_{L_{2}(\mathbb{R})}^{2},\quad\omega\in
H^{1}(\mathbb{R}),\ k\in\lbrack-\pi,\pi),\ \varepsilon>0.\nonumber
\end{align}
From (\ref{5.9}) and (\ref{1.16}) it follows that%
\begin{equation}
s(k,\varepsilon)[\omega,\omega]+\varepsilon^{2}\int_{\mathbb{R}}Q^{0}%
(x_{2})|\omega|^{2}dx_{2}\geq(c_{0}k^{2}+c_{3}\varepsilon^{2})\|\omega
\|_{L_{2}(\mathbb{R})}^{2},\quad\omega\in H^{1}(\mathbb{R}),\ k\in\lbrack
-\pi,\pi),\ \varepsilon>0, \label{5.10}%
\end{equation}
which implies (\ref{5.3}). $\blacksquare$

\begin{lemma}
\label{lem5.2}We have%
\begin{equation}
\|(D_{1}+k)(A(k,\varepsilon)+\varepsilon^{2}Q)^{-\frac{1}{2}}\|_{L_{2}%
(\Omega)\rightarrow L_{2}(\Omega)}\leq c_{0}^{-\frac{1}{2}},\quad k\in
\lbrack-\pi,\pi),\ \varepsilon>0, \label{5.11}%
\end{equation}%
\begin{equation}
\|\varepsilon D_{2}(A(k,\varepsilon)+\varepsilon^{2}Q)^{-\frac{1}{2}}%
\|_{L_{2}(\Omega)\rightarrow L_{2}(\Omega)}\leq c_{0}^{-\frac{1}{2}},\quad
k\in\lbrack-\pi,\pi),\ \varepsilon>0, \label{5.12}%
\end{equation}%
\begin{equation}
\|(D_{1}+k)A_{1}(k)^{-\frac{1}{2}}\|_{L_{2}(\Omega)\rightarrow L_{2}(\Omega
)}\leq c_{0}^{-\frac{1}{2}},\quad k\in\lbrack-\pi,\pi),\ k\neq0, \label{5.13}%
\end{equation}%
\begin{equation}
\|\varepsilon D_{2}(S(k,\varepsilon)+\varepsilon^{2}Q^{0})^{-\frac{1}{2}%
}P\|_{L_{2}(\Omega)\rightarrow L_{2}(\Omega)}\leq c_{0}^{-\frac{1}{2}},\quad
k\in\lbrack-\pi,\pi). \label{5.14}%
\end{equation}

\end{lemma}

\textbf{Proof.} By (\ref{1.1}) and (\ref{2.11}), we obtain%
\begin{align}
c_{0}\int_{\Omega}\left(  |(D_{1}+k)u|^{2}+\varepsilon^{2}|D_{2}u|^{2}\right)
\,d\mathbf{x}  &  \leq a(k,\varepsilon)[u,u]+\varepsilon^{2}\int_{\Omega
}Q(\mathbf{x})|u|^{2}\,d\mathbf{x}=\label{5.15}\\
&  =\|(A(k,\varepsilon)+\varepsilon^{2}Q)^{\frac{1}{2}}u\|_{L_{2}(\Omega)}%
^{2},\quad\ u\in\widetilde{H}^{1}(\Omega),\ k\in\lbrack-\pi,\pi),\ \varepsilon
>0.\nonumber
\end{align}
Substituting $u=(A(k,\varepsilon)+\varepsilon^{2}Q)^{-\frac{1}{2}}v,\,v\in
L_{2}(\Omega),$ in (\ref{5.15}), we derive (\ref{5.11}) and (\ref{5.12}).

Similarly, (\ref{5.6}) implies (\ref{5.13}). Finally, (\ref{5.14}) follows
from the estimate%
\begin{align*}
c_{0}\int_{\mathbb{R}}(k^{2}|\omega|^{2}+\varepsilon^{2}|D_{2}\omega
|^{2})\,dx_{2}  &  \leq s(k,\varepsilon)[\omega,\omega]+\varepsilon^{2}%
\int_{\mathbb{R}}Q^{0}(x_{2})|\omega|^{2}dx_{2}=\\
&  =\|(S(k,\varepsilon)+\varepsilon^{2}Q^{0})^{\frac{1}{2}}\omega
\|_{L_{2}(\mathbb{R})}^{2},\quad\omega\in H^{1}(\mathbb{R}),\ k\in\lbrack
-\pi,\pi),\ \varepsilon>0,
\end{align*}
see (\ref{5.9}). $\blacksquare$

\subsection{Estimates for generalized resolvents multiplied by projections}

Recall that $F(k)$ denotes the spectral projection of the operator $A_{1}(k)$
for the interval $[0,\delta]$, and $P$ denotes the orthogonal projection of
$L_{2}(\Omega)$ onto the subspace $\mathfrak{N}$ of functions depending only
on $x_{2}$.

\begin{lemma}
\label{lem5.3}Let $t=|k|\leq t_{0},\,\varepsilon>0.$ Then%
\begin{equation}
\|F(k)^{\perp}(A(k,\varepsilon)+\varepsilon^{2}Q)^{-\frac{1}{2}}%
\|_{L_{2}(\Omega)\rightarrow L_{2}(\Omega)}=\|(A(k,\varepsilon)+\varepsilon
^{2}Q)^{-\frac{1}{2}}F(k)^{\perp}\|_{L_{2}(\Omega)\rightarrow L_{2}(\Omega
)}\leq\delta^{-\frac{1}{2}}, \label{5.16}%
\end{equation}%
\begin{equation}
\|(A(k,\varepsilon)+\varepsilon^{2}Q)^{-\frac{1}{2}}P^{\perp}\|_{L_{2}%
(\Omega)\rightarrow L_{2}(\Omega)}\leq C_{3}:=\delta^{-\frac{1}{2}}+C_{1}%
c_{0}^{-\frac{1}{2}}, \label{5.17}%
\end{equation}%
\begin{equation}
\|F(k)^{\perp}(S(k,\varepsilon)+\varepsilon^{2}Q^{0})^{-\frac{1}{2}}%
P\|_{L_{2}(\Omega)\rightarrow L_{2}(\Omega)}\leq C_{4}:=C_{1}c_{0}^{-\frac
{1}{2}}. \label{5.18}%
\end{equation}

\end{lemma}

\textbf{Proof.} By (\ref{2.11}), (\ref{3.1}) and (\ref{1.7}), we have%
\[
a_{1}(k)[u,u]+\varepsilon^{2}c_{3}\|u\|_{L_{2}(\Omega)}^{2}\leq
a(k,\varepsilon)[u,u]+\varepsilon^{2}(Qu,u)_{L_{2}(\Omega)},\quad
u\in\widetilde{H}^{1}(\Omega),\ k\in\lbrack-\pi,\pi),\ \varepsilon>0.
\]
Hence, we get:%
\begin{equation}
\|(A_{1}(k)+\varepsilon^{2}c_{3}I)^{\frac{1}{2}}(A(k,\varepsilon
)+\varepsilon^{2}Q)^{-\frac{1}{2}}\|_{L_{2}(\Omega)\rightarrow L_{2}(\Omega
)}\leq1,\quad k\in\lbrack-\pi,\pi),\ \varepsilon>0. \label{5.19}%
\end{equation}
Therefore,%
\begin{align*}
\|F(k)^{\perp}(A(k,\varepsilon)+  &  \varepsilon^{2}Q)^{-\frac{1}{2}}%
\|_{L_{2}(\Omega)\rightarrow L_{2}(\Omega)}\leq\\
&  \leq\|F(k)^{\perp}(A_{1}(k)+\varepsilon^{2}c_{3}I)^{-\frac{1}{2}}%
\|_{L_{2}(\Omega)\rightarrow L_{2}(\Omega)}\|(A_{1}(k)+\varepsilon^{2}%
c_{3}I)^{\frac{1}{2}}(A(k,\varepsilon)+\varepsilon^{2}Q)^{-\frac{1}{2}}%
\|\leq\\
&  \leq\|F(k)^{\perp}(A_{1}(k)+\varepsilon^{2}c_{3}I)^{-\frac{1}{2}}%
\|_{L_{2}(\Omega)\rightarrow L_{2}(\Omega)}\leq\delta^{-\frac{1}{2}},
\end{align*}
which implies (\ref{5.16}).

To prove (\ref{5.17}), we use (\ref{3.28}), which gives $P^{\perp}%
=F(k)^{\perp}+t\Phi(k)$, and (\ref{3.32}). Then, by (\ref{5.16}) and
(\ref{5.1}), we get%
\begin{align*}
\| (A(k,\varepsilon)+\varepsilon^{2}Q)^{-\frac{1}{2}}P^{\perp}\|  &
_{L_{2}(\Omega)\rightarrow L_{2}(\Omega)} =\|(A(k,\varepsilon)+\varepsilon
^{2}Q)^{-\frac{1}{2}}(F(k)^{\perp}+t\Phi(k))\|\leq\\
&  \leq\delta^{-\frac{1}{2}}+C_{1}t(c_{0}t^{2}+c_{3}\varepsilon^{2}%
)^{-\frac{1}{2}}\leq\delta^{-\frac{1}{2}}+C_{1}c_{0}^{-\frac{1}{2}}%
=C_{3},\quad t=|k|\leq t_{0},\ \varepsilon>0.
\end{align*}

Finally, note that $P^{\perp}(S(k,\varepsilon)+\varepsilon^{2}Q^{0}%
)^{-\frac{1}{2}}P=0.$ Then, by (\ref{3.28}), $F(k)^{\perp}=P^{\perp}%
-t\Phi(k),$ and from (\ref{3.32}) and (\ref{5.3}) we obtain%
\begin{align*}
\|F(k)^{\perp}(S(k,\varepsilon)+\varepsilon^{2}Q^{0})^{-\frac{1}{2}}P\|  &
_{L_{2}(\Omega)\rightarrow L_{2}(\Omega)} =\|t\Phi(k)(S(k,\varepsilon
)+\varepsilon^{2}Q^{0})^{-\frac{1}{2}}P\|\leq\\
&  \leq C_{1}t(c_{0}t^{2}+c_{3}\varepsilon^{2})^{-\frac{1}{2}}\leq C_{1}%
c_{0}^{-\frac{1}{2}}=C_{4},\quad t=|k|\leq t_{0},\ \varepsilon>0.
\ \ \blacksquare
\end{align*}

Lemma \ref{lem5.3} together with Lemma \ref{lem5.1} imply the following corollary.

\begin{corollary}
\label{cor5.4}Let $t=|k|\leq t_{0},\,\varepsilon>0.$ Then%
\begin{equation}
\|(A(k,\varepsilon)+\varepsilon^{2}Q)^{-1}P^{\perp}\|_{L_{2}(\Omega
)\rightarrow L_{2}(\Omega)}\leq C_{3}c_{3}^{-\frac{1}{2}}\varepsilon^{-1},
\label{5.20}%
\end{equation}%
\begin{equation}
\|F(k)^{\perp}(S(k,\varepsilon)+\varepsilon^{2}Q^{0})^{-1}P\|_{L_{2}%
(\Omega)\rightarrow L_{2}(\Omega)}\leq C_{4}c_{3}^{-\frac{1}{2}}%
\varepsilon^{-1}. \label{5.21}%
\end{equation}

\end{corollary}

Indeed, (\ref{5.20}) follows from (\ref{5.1}) and (\ref{5.17}), while
(\ref{5.21}) is a consequence of (\ref{5.3}) and (\ref{5.18}).

\subsection{One more reduction}

Now we represent the operator under the norm sign in (\ref{4.5}) as%
\begin{align}
(A(k,\varepsilon)+\varepsilon^{2}Q)^{-1}  &  -(S(k,\varepsilon)+\varepsilon
^{2}Q^{0})^{-1}P =\label{5.22}\\
&  =(A(k,\varepsilon)+\varepsilon^{2}Q)^{-1}P^{\perp}-F(k)^{\perp
}(S(k,\varepsilon)+\varepsilon^{2}Q^{0})^{-1}P+\mathcal{J}(k,\varepsilon
),\nonumber
\end{align}
where
\begin{equation}
\mathcal{J}(k,\varepsilon):=(A(k,\varepsilon)+\varepsilon^{2}Q)^{-1}%
P-F(k)(S(k,\varepsilon)+\varepsilon^{2}Q^{0})^{-1}P. \label{5.23}%
\end{equation}
Estimates (\ref{5.20}), (\ref{5.21}) and identity (\ref{5.22}) show that
Theorem \ref{theor4.1} is a direct consequence of the following proposition.

\begin{proposition}
\label{prop5.5}Suppose that conditions of Theorem \textrm{{\ref{theor4.1}}}
are satisfied. Let $\mathcal{J}(k,\varepsilon)$ be the operator defined by
\textrm{{(\ref{5.23})}}. Then%
\begin{equation}
\|\mathcal{J}(k,\varepsilon)\|_{L_{2}(\Omega)\rightarrow L_{2}(\Omega)}%
\leq\widehat{C}\varepsilon^{-1},\quad|k|\leq t_{0},\ 0<\varepsilon\leq1.
\label{5.24}%
\end{equation}
The constant $\widehat{C}$ depends only on $c_{j},\,j=0,...,5.$
\end{proposition}

Clearly, the constant $\widetilde{C}$ in (\ref{4.5}) is given by
$\widetilde{C}= (C_{3} + C_{4})c_{3}^{-1/2} + \widehat{C}$. Proposition
\ref{prop5.5} (and Theorem \ref{theor4.1} with it) will be proved in Section
\ref{sect7}.

\section{Estimates for commutators\label{sect6}}

For the proof of Proposition \ref{prop5.5}, we need estimates for the
commutators of the operator $D_{2}$ with the generalized resolvent of
$A(k,\varepsilon)$ and with the projection $F(k)$.

\subsection{The commutator of $D_{2}$ and $(A(k,\varepsilon) + \varepsilon^{2}
Q)^{-1}$}

The following fact can be proved by analogy with the proof of Theorem 8.8 in
\cite{GT}.

\begin{lemma}
\label{lem6.1}Suppose that $g_{j}(\mathbf{x}),$ $j=1,2$, satisfy conditions of
the form \textrm{(\ref{1.1})} and \textrm{(\ref{1.3})} in $\Omega
=(0,1)\times\mathbb{R}.$ Suppose that $Q(\mathbf{x})$ satisfies conditions of
the form \textrm{(\ref{1.7})} and \textrm{(\ref{1.8})} in $\Omega$. Let
$u\in\widetilde{H}^{1}(\Omega)$ be a \textrm{(}weak\textrm{)} solution of the
equation%
\begin{equation}
(D_{1}+k)g_{1}(\mathbf{x})(D_{1}+k)u+\varepsilon^{2}D_{2}g_{2}(\mathbf{x}%
)D_{2}u+\varepsilon^{2}Q(\mathbf{x})u=f(\mathbf{x}), \label{6.1}%
\end{equation}
where $f\in L_{2}(\Omega)$, which means that $u$ satisfies the identity%
\begin{equation}
\int_{\Omega}\left(  g_{1}(\mathbf{x})(D_{1}+k)u\overline{(D_{1}%
+k)w}+\varepsilon^{2}g_{2}(\mathbf{x})D_{2}u\overline{D_{2}w}+\varepsilon
^{2}Q(\mathbf{x})u\overline{w}\right)  \,d\mathbf{x}=\int_{\Omega}f
\overline{w}\,d\mathbf{x},\ \forall w\in\widetilde{H}^{1}(\Omega). \label{6.2}%
\end{equation}
Then%
\begin{equation}
D_{2}u\in\widetilde{H}^{1}(\Omega). \label{6.3}%
\end{equation}

\end{lemma}

Let us prove the following lemma.

\begin{lemma}
\label{lem6.2}Let $|k|\leq t_{0},\,\varepsilon>0.$ The operator%
\begin{equation}
K_{1}(k,\varepsilon):=\operatorname*{clos}[(A(k,\varepsilon)+\varepsilon
^{2}Q)^{-1},D_{2}] \label{6.4}%
\end{equation}
is bounded in $L_{2}(\Omega)$, and satisfies the estimate%
\begin{equation}
\|K_{1}(k,\varepsilon)P^{\perp}\|_{L_{2}(\Omega)\rightarrow L_{2}(\Omega
)}=\|P^{\perp}K_{1}(k,\varepsilon)\|_{L_{2}(\Omega)\rightarrow L_{2}(\Omega
)}\leq C_{5}(c_{0}k^{2}+c_{3}\varepsilon^{2})^{-\frac{1}{2}}, \label{6.5}%
\end{equation}
where the constant $C_{5}$ is defined below in \textrm{{(\ref{6.9})}} and
depends only on $c_{0},c_{1},c_{2},c_{3},c_{5}.$
\end{lemma}

\textbf{Proof.} Let $u,v\in\widetilde{H}^{1}(\Omega)$ be such that
$D_{2}u,\ D_{2}v\in\widetilde{H}^{1}(\Omega)$. Then we have%
\begin{align}
a(k,\varepsilon)  &  [D_{2}u,v]+\varepsilon^{2}(QD_{2}u,v)_{L_{2}(\Omega
)}-a(k,\varepsilon)[u,D_{2}v]-\varepsilon^{2}(Qu,D_{2}v)_{L_{2}(\Omega)}
=\label{6.6}\\
&  =\int_{\Omega}\left(  g_{1}(\mathbf{x})((D_{1}+k)D_{2}u)\overline
{(D_{1}+k)v}+\varepsilon^{2}g_{2}(\mathbf{x})D_{2}^{2}u\overline{D_{2}%
v}\right)  \,d\mathbf{x}-\nonumber\\
&  \quad-\int_{\Omega}\left(  g_{1}(\mathbf{x})(D_{1}+k)u (\overline
{(D_{1}+k)D_{2}v})+\varepsilon^{2}g_{2}(\mathbf{x})D_{2}u\overline{D_{2}^{2}%
v}\right)  \,d\mathbf{x}+\varepsilon^{2}\int_{\Omega}Q(\mathbf{x})\left(
(D_{2}u)\overline{v}-u\overline{D_{2}v}\right)  \,d\mathbf{x} =\nonumber\\
&  =-\int_{\Omega}\left(  (D_{2}g_{1})(\mathbf{x})(D_{1}+k)u\overline
{(D_{1}+k)v}+\varepsilon^{2}(D_{2}g_{2})(\mathbf{x})D_{2}u\overline{D_{2}%
v}\right)  \,d\mathbf{x}-\varepsilon^{2}\int_{\Omega}(D_{2}Q)(\mathbf{x}%
)u\overline{v}\,d\mathbf{x}.\nonumber
\end{align}
Here we have integrated by parts in direction $x_{2}$. Now we put%
\begin{equation}
u=(A(k,\varepsilon)+\varepsilon^{2}Q)^{-1}\varphi,\ \ \ v=(A(k,\varepsilon
)+\varepsilon^{2}Q)^{-1}\psi,\ \ \ \ \varphi, \psi\in L_{2}(\Omega).
\label{6.7}%
\end{equation}
By Lemma \ref{lem6.1}, we have $u,v\in\widetilde{H}^{1}(\Omega),\,D_{2}%
u,D_{2}v\in\widetilde{H}^{1}(\Omega).$

Then from (\ref{6.6}) we obtain%
\begin{align*}
(D_{2}(A(k,  &  \varepsilon)+\varepsilon^{2}Q)^{-1}\varphi,\psi)_{L_{2}%
(\Omega)}-(\varphi,D_{2}(A(k,\varepsilon)+\varepsilon^{2}Q)^{-1}\psi
)_{L_{2}(\Omega)} =\\
&  =-\int_{\Omega}(D_{2}g_{1})(\mathbf{x})\left(  (D_{1}+k)(A(k,\varepsilon
)+\varepsilon^{2}Q)^{-1}\varphi\right)  \left(  \overline{(D_{1}%
+k)(A(k,\varepsilon)+\varepsilon^{2}Q)^{-1}\psi}\right)  \,d\mathbf{x}-\\
&  \quad-\int_{\Omega}\varepsilon^{2}(D_{2}g_{2})(\mathbf{x})\left(
D_{2}(A(k,\varepsilon)+\varepsilon^{2}Q)^{-1}\varphi\right)  \left(
\overline{D_{2}(A(k,\varepsilon)+\varepsilon^{2}Q)^{-1}\psi}\right)
\,d\mathbf{x}-\\
&  \quad-\int_{\Omega}\varepsilon^{2}(D_{2}Q)(\mathbf{x})\left(
(A(k,\varepsilon)+\varepsilon^{2}Q)^{-1}\varphi\right)  \left(  \overline
{(A(k,\varepsilon)+\varepsilon^{2}Q)^{-1}\psi}\right)  \,d\mathbf{x}.
\end{align*}
It follows that%
\begin{align}
K_{1}(k,\varepsilon)  &  =-((D_{1}+k)(A(k,\varepsilon)+\varepsilon^{2}%
Q)^{-1})^{\ast}[D_{2}g_{1}]\left(  (D_{1}+k)(A(k,\varepsilon)+\varepsilon
^{2}Q)^{-1}\right)  -\label{6.8}\\
&  \quad-(\varepsilon D_{2}(A(k,\varepsilon)+\varepsilon^{2}Q)^{-1})^{\ast
}[D_{2}g_{2}](\varepsilon D_{2}(A(k,\varepsilon)+\varepsilon^{2}%
Q)^{-1})-\nonumber\\
&  \quad-\varepsilon^{2}(A(k,\varepsilon)+\varepsilon^{2}Q)^{-1}%
[D_{2}Q](A(k,\varepsilon)+\varepsilon^{2}Q)^{-1}.\nonumber
\end{align}
By (\ref{1.3}), (\ref{1.8}), (\ref{5.1}), (\ref{5.11}) and (\ref{5.12}), the
operator on the right is bounded in $L_{2}(\Omega)$.

Next, from (\ref{1.3}), (\ref{1.8}), (\ref{5.1}), (\ref{5.11}), (\ref{5.12}),
(\ref{5.17}) and (\ref{6.8}) it follows that%
\begin{align*}
\|K_{1}  &  (k,\varepsilon)P^{\perp}\|_{L_{2}(\Omega)\rightarrow L_{2}%
(\Omega)}\leq\\
&  \leq\|(D_{1}+k)(A(k,\varepsilon)+\varepsilon^{2}Q)^{-\frac{1}{2}}%
\|^{2}_{L_{2}(\Omega)\rightarrow L_{2}(\Omega)} \|(A(k,\varepsilon
)+\varepsilon^{2}Q)^{-\frac{1}{2}}\|_{L_{2}(\Omega)\rightarrow L_{2}(\Omega
)}\times\\
&  \qquad\times\|D_{2}g_{1}\|_{L_{\infty}}\|(A(k,\varepsilon)+\varepsilon
^{2}Q)^{-\frac{1}{2}}P^{\perp}\|_{L_{2}(\Omega)\rightarrow L_{2}(\Omega)}+\\
&  \ \ +\|\varepsilon D_{2}(A(k,\varepsilon)+\varepsilon^{2}Q)^{-\frac{1}{2}%
}\|_{L_{2}(\Omega)\rightarrow L_{2}(\Omega)}^{2}\|(A(k,\varepsilon
)+\varepsilon^{2}Q)^{-\frac{1}{2}}\|_{L_{2}(\Omega)\rightarrow L_{2}(\Omega
)}\times\\
&  \qquad\times\|D_{2}g_{2}\|_{L_{\infty}}\|(A(k,\varepsilon)+\varepsilon
^{2}Q)^{-\frac{1}{2}}P^{\perp}\|_{L_{2}(\Omega)\rightarrow L_{2}(\Omega)}+\\
&  \ \ +\varepsilon^{2}\|(A(k,\varepsilon)+\varepsilon^{2}Q)^{-1}%
\|_{L_{2}(\Omega)\rightarrow L_{2}(\Omega)}^{\frac{3}{2}}\|D_{2}%
Q\|_{L_{\infty}}\|(A(k,\varepsilon)+\varepsilon^{2}Q)^{-\frac{1}{2}}P^{\perp
}\|_{L_{2}(\Omega)\rightarrow L_{2}(\Omega)}\leq\\
&  \leq2c_{0}^{-1}(c_{0}k^{2}+c_{3}\varepsilon^{2})^{-\frac{1}{2}}c_{2}%
C_{3}+\varepsilon^{2}(c_{0}k^{2}+c_{3}\varepsilon^{2})^{-\frac{3}{2}}%
c_{5}C_{3},
\end{align*}
which implies (\ref{6.5}) with%
\begin{equation}
C_{5}=2c_{0}^{-1}c_{2}C_{3}+c_{3}^{-1}c_{5}C_{3}. \ \ \blacksquare\label{6.9}%
\end{equation}

\subsection{The commutator of $D_{2}$ and $F(k)$}

In order to calculate the commutator $[D_{2},F(k)]$, we represent the spectral
projection $F(k)$ as the integral of the resolvent $(A_{1}(k)-zI)^{-1}$ over
the contour $\Gamma_{\delta}$ that envelopes the real interval $[0,\delta]$
equidistantly with distance $\delta$. Recall that, for $|k|\leq t_{0}$, the
interval $(\delta,3\delta)$ is free of the spectrum of $A_{1}(k)$; therefore,%
\begin{equation}
\|(A_{1}(k)-zI)^{-1}\|_{L_{2}(\Omega)\rightarrow L_{2}(\Omega)}\leq\delta
^{-1},\,\ \ |k|\leq t_{0},\quad z\in\Gamma_{\delta}. \label{6.10}%
\end{equation}
We have the following representation%
\begin{equation}
F(k)= - \frac{1}{2\pi i}\int_{\Gamma_{\delta}}(A_{1}(k)-zI)^{-1}dz.
\label{6.11}%
\end{equation}

\begin{lemma}
\label{lem6.3}Let $t=|k|\leq t_{0}$ and $\varepsilon>0.$ Then the operator%
\begin{equation}
K_{2}(k):=\operatorname*{clos}[D_{2},F(k)] \label{6.12}%
\end{equation}
is bounded in $L_{2}(\Omega)$, and%
\begin{equation}
\|K_{2}(k)\|_{L_{2}(\Omega)\rightarrow L_{2}(\Omega)}\leq C_{6}, \label{6.13}%
\end{equation}%
\begin{equation}
\|K_{2}(k)F(k)\|_{L_{2}(\Omega)\rightarrow L_{2}(\Omega)}\leq C_{7}t,
\label{6.14}%
\end{equation}
where the constants $C_{6}$ and $C_{7}$ are defined below in
\textrm{{(\ref{6.21})}} and \textrm{{(\ref{6.27})}} respectively. These
constants depend only on $c_{0},c_{1},c_{2}.$
\end{lemma}

\textbf{Proof.} We start with calculating the operator%
\begin{equation}
K_{3}(k,z):=\operatorname*{clos}[D_{2},(A_{1}(k)-zI)^{-1}]. \label{6.15}%
\end{equation}
Suppose that functions $u,v\in\operatorname*{Dom}a_{1}(k)$ are such that
$D_{2}u,D_{2}v\in\operatorname*{Dom}a_{1}(k)$. Then%
\begin{align}
a_{1}(k)  &  [D_{2}u,v]-z(D_{2}u,v)_{L_{2}(\Omega)}-a_{1}(k)[u,D_{2}%
v]+z(u,D_{2}v)_{L_{2}(\Omega)}=\label{6.16}\\
&  =\int_{\Omega}\left(  g_{1}(\mathbf{x})((D_{1}+k)D_{2}u)\overline
{(D_{1}+k)v}-g_{1}(\mathbf{x})((D_{1}+k)u) \overline{(D_{1}+k)D_{2}v}\right)
\,d\mathbf{x}+\nonumber\\
&  \qquad+z\int_{\Omega}\left(  u\overline{D_{2}v}- (D_{2}u) \overline
{v}\right)  \,d\mathbf{x}=\nonumber\\
&  =-\int_{\Omega}(D_{2}g_{1})(\mathbf{x})((D_{1}+k)u)\overline{(D_{1}%
+k)v}\,d\mathbf{x}.\nonumber
\end{align}
We have integrated by parts in direction $x_{2}$. We put $u=(A_{1}%
(k)-zI)^{-1}\varphi,\,v=(A_{1}(k)-zI)^{-1}\psi,$ where $\varphi,\psi
\in\widetilde{C}_{0}^{\infty}(\Omega)$. Here $\widetilde{C}_{0}^{\infty
}(\Omega)$ is a class of infinitely smooth functions in $\Omega$ periodic in
$x_{1}$ (which means that 1-periodic extension of them to $\mathbb{R}^{2}$ is
$C^{\infty}$ in $\mathbb{R}^{2}$) and equal to zero for large $|x_{2}|.$ By
using (\ref{1.3}), it is easily seen that $u$ and $v$ satisfy the required
conditions. Then from (\ref{6.16}) we obtain for all $\varphi,\psi
\in\widetilde{C}_{0}^{\infty}(\Omega)$:%
\begin{align*}
(D_{2}(A_{1}(k)  &  -zI)^{-1}\varphi,\psi)_{L_{2}(\Omega)}-(\varphi
,D_{2}(A_{1}(k)-zI)^{-1}\psi)_{L_{2}(\Omega)}=\\
&  =-\int_{\Omega}(D_{2}g_{1})(\mathbf{x})\left(  (D_{1}+k)(A_{1}%
(k)-zI)^{-1}\varphi\right)  \left(  \overline{(D_{1}+k)(A_{1}(k)-zI)^{-1}\psi
}\right)  \,d\mathbf{x}.
\end{align*}
Since $\widetilde{C}_{0}^{\infty}(\Omega)$ is dense in $L_{2}(\Omega)$, it
follows that%
\begin{equation}
K_{3}(k,z)=-((D_{1}+k)(A_{1}(k)-zI)^{-1})^{\ast}[D_{2}g_{1}]((D_{1}%
+k)(A_{1}(k)-zI)^{-1}). \label{6.17}%
\end{equation}
Then, by (\ref{1.3}) and (\ref{5.13}), we have%
\begin{align*}
\|K_{3}(k,z)\|_{L_{2}(\Omega)\rightarrow L_{2}(\Omega)}  &  \leq\|D_{2}%
g_{1}\|_{L_{\infty}}\|(D_{1}+k)(A_{1}(k)-zI)^{-1}\|_{L_{2}(\Omega)\rightarrow
L_{2}(\Omega)}^{2}\leq\\
&  \leq c_{2}\|(D_{1}+k)A_{1}(k)^{-\frac{1}{2}}\|_{L_{2}(\Omega)\rightarrow
L_{2}(\Omega)}^{2}\|A_{1}(k)^{\frac{1}{2}}(A_{1}(k)-zI)^{-1}\|_{L_{2}%
(\Omega)\rightarrow L_{2}(\Omega)}^{2}\leq\\
&  \leq c_{2}c_{0}^{-1}\|A_{1}(k)^{\frac{1}{2}}(A_{1}(k)-zI)^{-1}%
\|_{L_{2}(\Omega)\rightarrow L_{2}(\Omega)}^{2}.
\end{align*}
Next, for $|k|\leq t_{0}$ we have:%
\begin{equation}
\|A_{1}(k)^{\frac{1}{2}}(A_{1}(k)-zI)^{-1}\|_{L_{2}(\Omega)\rightarrow
L_{2}(\Omega)}\leq\max\{{\lambda^{\frac{1}{2}}|\lambda-z|^{-1}:\lambda
\in\lbrack0,\delta]\cup\lbrack3\delta,\infty),\,z\in\Gamma_{\delta}%
\}}=:c(\delta). \label{6.18}%
\end{equation}
Hence,%
\begin{equation}
\|K_{3}(k,z)\|_{L_{2}(\Omega)\rightarrow L_{2}(\Omega)}\leq c_{2}c_{0}%
^{-1}c(\delta)^{2},\quad z\in\Gamma_{\delta},\ |k|\leq t_{0}. \label{6.19}%
\end{equation}

Now, from (\ref{6.11}), (\ref{6.12}) and (\ref{6.15}) it follows that%
\begin{equation}
K_{2}(k)=- \frac{1}{2\pi i}\int_{\Gamma_{\delta}}K_{3}(k,z)dz. \label{6.20}%
\end{equation}
Combining (\ref{6.19}) with (\ref{6.20}), we obtain%
\begin{equation}
\|K_{2}(k)\|_{L_{2}(\Omega)\rightarrow L_{2}(\Omega)}\leq(2\delta+2\pi
\delta)(2\pi)^{-1}c_{2}c_{0}^{-1}c(\delta)^{2}=:C_{6},\quad|k|\leq t_{0},
\label{6.21}%
\end{equation}
which proves (\ref{6.13}).

Next, by (\ref{6.17}), (\ref{1.3}), (\ref{5.13}), (\ref{6.18}), (\ref{6.10}),
we have%
\begin{align}
\|K_{3}(k,z)F(k)\|  &  _{L_{2}(\Omega)\rightarrow L_{2}(\Omega)}%
\leq\label{6.22}\\
&  \leq\|D_{2}g_{1}\|_{L_{\infty}}\|(D_{1}+k)(A_{1}(k)-zI)^{-1}\|_{L_{2}%
(\Omega)\rightarrow L_{2}(\Omega)}\times\nonumber\\
&  \qquad\times\|(D_{1}+k)(A_{1}(k)-zI)^{-1}F(k)\|_{L_{2}(\Omega)\rightarrow
L_{2}(\Omega)}\leq\nonumber\\
&  \leq c_{2}\|(D_{1}+k)A_{1}(k)^{-\frac{1}{2}}\|_{L_{2}(\Omega)\rightarrow
L_{2}(\Omega)}^{2}\|A_{1}(k)^{\frac{1}{2}}(A_{1}(k)-zI)^{-1}\|_{L_{2}%
(\Omega)\rightarrow L_{2}(\Omega)}\times\nonumber\\
&  \qquad\times\|(A_{1}(k)-zI)^{-1}\|_{L_{2}(\Omega)\rightarrow L_{2}(\Omega
)}\|A_{1}(k)^{\frac{1}{2}}F\left(  k\right)  \|_{L_{2}(\Omega)\rightarrow
L_{2}(\Omega)}\nonumber\\
&  \leq c_{2}c_{0}^{-1}c(\delta)\delta^{-1}\|A_{1}(k)^{\frac{1}{2}%
}F(k)\|_{L_{2}(\Omega)\rightarrow L_{2}(\Omega)},\quad|k|\leq t_{0}%
,\ z\in\Gamma_{\delta}.\nonumber
\end{align}
Recalling (\ref{3.30}) and (\ref{3.33}), we obtain%
\begin{align}
\|A_{1}(k)^{\frac{1}{2}}F(k)u\|_{L_{2}(\Omega)}^{2}  &  =(A_{1}%
(k)F(k)u,u)_{L_{2}(\Omega)}=((t^{2}SP+t^{3}\Psi(k))u,u)_{L_{2}(\Omega)}%
\leq\label{6.23}\\
&  \leq(t^{2}\|S\|_{\mathfrak{N}\rightarrow\mathfrak{N}}+C_{2}t^{3}%
)\|u\|_{L_{2}(\Omega)}^{2},\quad t=|k|\leq t_{0}.\nonumber
\end{align}
Since $S$ is the operator of multiplication by the function $g_{1}^{0}(x_{2}%
)$, by (\ref{1.13}), we have
\begin{equation}
\|S\|_{\mathfrak{N}\rightarrow\mathfrak{N}}\leq c_{1}. \label{6.24}%
\end{equation}
From (\ref{6.23}) and (\ref{6.24}) it follows that%
\begin{equation}
\|A_{1}(k)^{\frac{1}{2}}F(k)\|_{L_{2}(\Omega)\rightarrow L_{2}(\Omega)}\leq
t(c_{1}+C_{2}t_{0})^{\frac{1}{2}},\,\ \ t=|k|\leq t_{0}. \label{6.25}%
\end{equation}
Combining (\ref{6.22}) and (\ref{6.25}), we arrive at the estimate
\begin{equation}
\|K_{3}(k,z)F(k)\|_{L_{2}(\Omega)\rightarrow L_{2}(\Omega)}\leq c_{2}%
c_{0}^{-1}c(\delta)\delta^{-1}(c_{1}+C_{2}t_{0})^{\frac{1}{2}}%
t,\,\ \ t=|k|\leq t_{0},\,z\in\Gamma_{\delta}. \label{6.26}%
\end{equation}
Finally, from (\ref{6.26}) and (\ref{6.20}) we obtain (\ref{6.14}) with the
constant%
\begin{equation}
C_{7}=(2\delta+2\pi\delta)(2\pi)^{-1}c_{2}c_{0}^{-1}c(\delta)\delta^{-1}%
(c_{1}+C_{2}t_{0})^{\frac{1}{2}}.\ \ \blacksquare\label{6.27}%
\end{equation}

\bigskip

In the sequel, we shall need the following statement.

\begin{lemma}
\label{lem6.4}We have%
\begin{equation}
F(k)\widetilde{H}^{1}(\Omega)\subset\widetilde{H}^{1}(\Omega). \label{6.28}%
\end{equation}

\end{lemma}

\textbf{Proof.} Let $u\in\widetilde{H}^{1}(\Omega)$. Then, obviously,
$u\in\operatorname*{Dom}a_{1}(k)=\int_{\mathbb{R}}\oplus\widetilde{H}%
^{1}(0,1)dx_{2}.$

Since $F(k)\operatorname*{Dom}a_{1}(k)\subset\operatorname*{Dom}a_{1}(k)$,
then $F(k)u\in\operatorname*{Dom}a_{1}(k).$ This means that $F(k)u\in
L_{2}(\Omega)$, $D_{1}(F(k)u)\in L_{2}(\Omega)$ and $F(k)u$ satisfies periodic
boundary conditions on $\partial\Omega$. Next,%
\[
D_{2}(F(k)u)=F(k)D_{2}u+[D_{2},F(k)]u\in L_{2}(\Omega),
\]
because $D_{2}u\in L_{2}(\Omega)$ and the operators $F(k)$ and $[D_{2},F(k)]$
are bounded in $L_{2}(\Omega)$ (see Lemma \ref{lem6.3}). It follows that
$F(k)u\in\widetilde{H}^{1}(\Omega).$ $\blacksquare$

\section{Proof of Proposition \ref{prop5.5}\label{sect7}}

\subsection{Abstract lemma}

For the study of the operator (\ref{5.23}), we wish to apply the identity%
\begin{equation}
\mathcal{J}(k,\varepsilon)=(A(k,\varepsilon)+\varepsilon^{2}Q)^{-1}%
(S(k,\varepsilon)P+\varepsilon^{2}Q^{0}P-A(k,\varepsilon)F(k)-\varepsilon
^{2}QF(k))(S(k,\varepsilon)+\varepsilon^{2}Q^{0})^{-1}P. \label{7.1}%
\end{equation}
However, the operator in the central brackets makes no sense. In order to
obtain a precise version of (\ref{7.1}), we need the following statement of
abstract nature (see Lemma 7.1 in \cite{Su}).

\begin{lemma}
\label{lem7.1}Let $\mathfrak{t}_{1}$ and $\mathfrak{t}_{2}$ be two closed and
densely defined sesquilinear forms in a separable Hilbert space $\mathfrak{H}%
$. Suppose that both quadratic forms $\mathfrak{t}_{1}[u,u]$ and
$\mathfrak{t}_{2}[u,u]$ are positive definite. Let $T_{1}$ and $T_{2}$ be the
selfadjoint operators in $\mathfrak{H}$ corresponding to the forms
$\mathfrak{t}_{1}$ and $\mathfrak{t}_{2}$ respectively. Let $\mathcal{P}$ and
$\mathcal{Q}$ be some orthogonal projections in $\mathfrak{H}$. Let the
following conditions be satisfied:

%\begin{enumerate}
%\item
$1^{\circ}$. $\mathfrak{d}:=\operatorname*{Dom}\mathfrak{t}_{1}\subset
\operatorname*{Dom}\mathfrak{t}_{2}$;

$2^{\circ}$. $\mathcal{P}\operatorname*{Dom}\mathfrak{t}_{2}\subset
\operatorname*{Dom}\mathfrak{t}_{1}$;

$3^{\circ}$. $\mathcal{Q}\mathfrak{d}\subset\mathfrak{d}$;

$4^{\circ}$. $\mathcal{P}$ commutes with $T_{2}$;

$5^{\circ}$. the following representation holds%
\begin{equation}
\mathfrak{t}_{2}[\mathcal{P}u,v]-\mathfrak{t}_{1}[\mathcal{Q}u,v]=(G_{0}%
u,Gv)_{\mathfrak{G}},\,\ \ u,v\in\mathfrak{d}, \label{7.2}%
\end{equation}
%\end{enumerate}
where $G_{0},G:\mathfrak{H}\rightarrow\mathfrak{G}$ are densely defined linear
operators acting from $\mathfrak{H}$ to some separable Hilbert space
$\mathfrak{G}$. Moreover, $\mathfrak{d}\subset\operatorname*{Dom}%
G_{0},\,\mathfrak{d}\subset\operatorname*{Dom}G$, and the operators
\[
GT_{1}^{-1}:\mathfrak{H}\rightarrow\mathfrak{G},\,\ \ \ \ G_{0}T_{2}%
^{-1}\mathcal{P}:\mathfrak{H}\rightarrow\mathfrak{G}%
\]
are bounded. Then we have%
\begin{equation}
T_{1}^{-1}\mathcal{P}-\mathcal{Q}T_{2}^{-1}\mathcal{P}=(GT_{1}^{-1})^{\ast
}(G_{0}T_{2}^{-1}\mathcal{P}). \label{7.3}%
\end{equation}

\end{lemma}

\bigskip

Let $|k|\leq t_{0},\,\varepsilon>0$. We shall apply Lemma \ref{lem7.1} with
\begin{align*}
\mathfrak{H}  &  =L_{2}(\Omega),\ \ \mathcal{P}=P,\ \ \mathcal{Q}=F(k),\\
\mathfrak{t}_{1}[u,v]  &  =a(k,\varepsilon)[u,v]+\varepsilon^{2}%
(Qu,v)_{L_{2}(\Omega)},\\
\mathfrak{t}_{2}[u,v]  &  =s(k,\varepsilon)[Pu,Pv]+\varepsilon^{2}%
(Q^{0}u,v)_{L_{2}(\Omega)},\\
\mathfrak{d}  &  =\widetilde{H}^{1}(\Omega)\text{ \ and \ }\operatorname*{Dom}%
\mathfrak{t}_{2}=\{{u\in L_{2}(\Omega):Pu\in H^{1}(\mathbb{R})\}}.
\end{align*}
Obviously, conditions $1^{\circ}$, $2^{\circ}$ and $4^{\circ}$ are satisfied.
Condition $3^{\circ}$ is also satisfied, by Lemma \ref{lem6.4}.

Let us check condition $5^{\circ}$. The left-hand side of (\ref{7.2}) takes
the form
\begin{align}
\mathcal{I}[u,v]  &  :=s(k,\varepsilon)[Pu,Pv]+\varepsilon^{2}(Q^{0}%
Pu,v)_{L_{2}(\Omega)}-a(k,\varepsilon)[F(k)u,v]-\varepsilon^{2}%
(QF(k)u,v)_{L_{2}(\Omega)}=\label{7.4}\\
&  =\mathcal{I}_{1}[u,v]+\mathcal{I}_{2}[u,v]+\mathcal{I}_{3}[u,v],\nonumber
\end{align}
where, in accordance with (\ref{2.11}), (\ref{3.1}) and (\ref{4.1}),%
\begin{align}
\mathcal{I}_{1}[u,v]  &  =\varepsilon^{2}(Q^{0}Pu-QF(k)u,v)_{L_{2}(\Omega
)},\label{7.5}\\
\mathcal{I}_{2}[u,v]  &  =t^{2}(SPu,Pv)_{L_{2}(\Omega)}-a_{1}%
(k)[F(k)u,v],\label{7.6}\\
\mathcal{I}_{3}[u,v]  &  =(g_{2}^{0}\varepsilon D_{2}Pu,\varepsilon
D_{2}Pv)_{L_{2}(\Omega)}-(g_{2}\varepsilon D_{2}F(k)u,\varepsilon
D_{2}v)_{L_{2}(\Omega)}. \label{7.7}%
\end{align}

\subsection{The form $\mathcal{I}_{1}$}

We represent the form (\ref{7.5}) as%
\begin{equation}
\mathcal{I}_{1}[u,v]=\mathcal{I}_{1}^{\prime}[u,v]+\mathcal{I}_{1}%
^{\prime\prime}[u,v], \label{7.8}%
\end{equation}
where%
\begin{align}
\mathcal{I}_{1}^{\prime}[u,v]  &  =\varepsilon^{2}((Q^{0}-Q)Pu,v)_{L_{2}%
(\Omega)},\label{7.9}\\
\mathcal{I}_{1}^{\prime\prime}[u,v]  &  =\varepsilon^{2}(Q(P-F(k))u,v)_{L_{2}%
(\Omega)}. \label{7.10}%
\end{align}
We have
\begin{align}
\mathcal{I}_{1}^{\prime}[u,v]  &  =\varepsilon^{2}(Q^{0}Pu,Pv)_{L_{2}(\Omega
)}-\varepsilon^{2}(QPu,Pv)_{L_{2}(\Omega)}-\varepsilon^{2}(QPu,P^{\perp
}v)_{L_{2}(\Omega)}=\label{7.11}\\
&  =-\varepsilon^{2}(QPu,P^{\perp}v)_{L_{2}(\Omega)},\nonumber
\end{align}
since $PQP=Q^{0}P$. Then%
\begin{equation}
\mathcal{I}_{1}^{\prime}[u,v]=(G_{01}^{\prime}u,G_{1}^{\prime}v)_{L_{2}%
(\Omega)}, \label{7.12}%
\end{equation}
where $G_{01}^{\prime}=-\varepsilon^{\frac{3}{2}}QP,\,G_{1}^{\prime
}=\varepsilon^{\frac{1}{2}}P^{\perp}$. By (\ref{5.3}) and (\ref{1.7}), we
obtain:%
\begin{equation}
\|G_{01}^{\prime}(S(k,\varepsilon)P+\varepsilon^{2}Q^{0})^{-1}P\|_{L_{2}%
(\Omega)\rightarrow L_{2}(\Omega)}\leq\varepsilon^{\frac{3}{2}}%
\|Q\|_{L_{\infty}}(c_{0}k^{2}+c_{3}\varepsilon^{2})^{-1}\leq c_{4}c_{3}%
^{-1}\varepsilon^{-\frac{1}{2}}. \label{7.13}%
\end{equation}
Using (\ref{5.20}), we obtain that%
\begin{align}
\|G_{1}^{\prime}(A(k,\varepsilon)+\varepsilon^{2}Q)^{-1}\|_{L_{2}%
(\Omega)\rightarrow L_{2}(\Omega)}  &  \leq\varepsilon^{\frac{1}{2}}%
\|P^{\perp}(A(k,\varepsilon)+\varepsilon^{2}Q)^{-1}\|_{L_{2}(\Omega
)\rightarrow L_{2}(\Omega)}\leq\label{7.14}\\
&  \leq C_{3}c_{3}^{-\frac{1}{2}}\varepsilon^{-\frac{1}{2}}.\nonumber
\end{align}

By (\ref{3.28}), the form (\ref{7.10}) can be written as%
\begin{equation}
\mathcal{I}_{1}^{\prime\prime}[u,v] = -\varepsilon^{2}(Qt\Phi(k)u,v)_{L_{2}%
(\Omega)}=(G_{01}^{\prime\prime}u,G_{1}^{\prime\prime}v)_{L_{2}(\Omega)},
\label{7.15}%
\end{equation}
where $G_{01}^{\prime\prime}=-\varepsilon t^{\frac{1}{2}}Q\Phi(k),$
$G_{1}^{\prime\prime}=\varepsilon t^{\frac{1}{2}}I.$

Using (\ref{3.32}), (\ref{5.3}) and (\ref{1.7}), we have%
\begin{align}
\|G_{01}^{\prime\prime}(S(k,\varepsilon)P+\varepsilon^{2}Q^{0})^{-1}%
P\|_{L_{2}(\Omega)\rightarrow L_{2}(\Omega)}  &  \leq\varepsilon t^{\frac
{1}{2}}\|Q\|_{L_{\infty}}C_{1}(c_{0}k^{2}+c_{3}\varepsilon^{2})^{-1}%
\leq\label{7.16}\\
&  \leq C_{1}c_{4}c_{0}^{-\frac{1}{4}}c_{3}^{-\frac{3}{4}}\varepsilon
^{-\frac{1}{2}}.\nonumber
\end{align}
Next, by (\ref{5.1}),
\begin{equation}
\|G_{1}^{\prime\prime}(A(k,\varepsilon)+\varepsilon^{2}Q)^{-1}\|_{L_{2}%
(\Omega)\rightarrow L_{2}(\Omega)}\leq\varepsilon t^{\frac{1}{2}}(c_{0}%
k^{2}+c_{3}\varepsilon^{2})^{-1}\leq c_{0}^{-\frac{1}{4}}c_{3}^{-\frac{3}{4}%
}\varepsilon^{-\frac{1}{2}}. \label{7.17}%
\end{equation}

\subsection{The form $\mathcal{I}_{2}$}

As for the form (\ref{7.6}), we apply (\ref{3.30}):%
\begin{align}
\mathcal{I}_{2}[u,v]  &  =((t^{2}SP-A_{1}(k)F(k))u,v)_{L_{2}(\Omega)}%
=-t^{3}(\Psi(k)u,v)_{L_{2}(\Omega)}=\label{7.18}\\
&  =(G_{02}u,G_{2}v)_{L_{2}(\Omega)},\nonumber
\end{align}
where $G_{02}=-t^{\frac{3}{2}}\Psi(k),$ $G_{2}=t^{\frac{3}{2}}I.$ Then, by
(\ref{3.33}) and (\ref{5.3}) we have%
\begin{align}
\|G_{02}(S(k,\varepsilon)P+\varepsilon^{2}Q^{0})^{-1}P\|_{L_{2}(\Omega
)\rightarrow L_{2}(\Omega)}  &  \leq t^{\frac{3}{2}}C_{2}(c_{0}k^{2}%
+c_{3}\varepsilon^{2})^{-1}\leq\label{7.19}\\
&  \leq C_{2}c_{0}^{-\frac{3}{4}}c_{3}^{-\frac{1}{4}}\varepsilon^{-\frac{1}%
{2}}.\nonumber
\end{align}
By (\ref{5.1}),%
\begin{equation}
\|G_{2}(A(k,\varepsilon)+\varepsilon^{2}Q)^{-1}\|_{L_{2}(\Omega)\rightarrow
L_{2}(\Omega)}\leq t^{\frac{3}{2}}(c_{0}k^{2}+c_{3}\varepsilon^{2})^{-1}\leq
c_{0}^{-\frac{3}{4}}c_{3}^{-\frac{1}{4}}\varepsilon^{-\frac{1}{2}}.
\label{7.20}%
\end{equation}

\subsection{The form $\mathcal{I}_{3}$}

We rewrite the form (\ref{7.7}) as%
\begin{equation}
\mathcal{I}_{3}[u,v]=\mathcal{I}_{3}^{\prime}[u,v]+\mathcal{I}_{4}[u,v],
\label{7.21}%
\end{equation}
where%
\begin{align}
\mathcal{I}_{3}^{\prime}[u,v]  &  =-(g_{2}\varepsilon D_{2}%
(F(k)-P)u,\varepsilon D_{2}v)_{L_{2}(\Omega)},\label{7.22}\\
\mathcal{I}_{4}[u,v]  &  =(g_{2}^{0}\varepsilon D_{2}Pu,\varepsilon
D_{2}Pv)_{L_{2}(\Omega)}-(g_{2}\varepsilon D_{2}Pu,\varepsilon D_{2}%
v)_{L_{2}(\Omega)}. \label{7.23}%
\end{align}
For the form (\ref{7.22}) we have%
\begin{equation}
\mathcal{I}_{3}^{\prime}[u,v]=(G_{03}u,G_{3}v)_{L_{2}(\Omega)}, \label{7.24}%
\end{equation}
where $G_{03}=-g_{2}\varepsilon^{\frac{1}{2}}D_{2}(F(k)-P),\,G_{3}%
=\varepsilon^{\frac{3}{2}}D_{2}.$

We represent $G_{03}$ as%
\begin{equation}
G_{03}=\widetilde{G}_{03}+\widehat{G}_{03}+\overset{\vee}{G}_{03},
\label{7.25}%
\end{equation}
where%
\begin{align*}
\widetilde{G}_{03}  &  =-g_{2}\varepsilon^{\frac{1}{2}}[D_{2},F(k)]F(k),\\
\widehat{G}_{03}  &  =-g_{2}\varepsilon^{\frac{1}{2}}[D_{2},F(k)]F(k)^{\perp
},\\
\overset{\vee}{G}_{03}  &  =-g_{2}\varepsilon^{\frac{1}{2}}(F(k)-P)D_{2}%
=-g_{2}\varepsilon^{\frac{1}{2}}t\Phi(k)D_{2}.
\end{align*}
Here we have taken (\ref{3.28}) into account. From (\ref{1.1}), (\ref{6.14})
and (\ref{5.3}) it follows that%
\begin{align}
\|\widetilde{G}_{03}  &  (S(k,\varepsilon)P+\varepsilon^{2}Q^{0}%
)^{-1}P\|_{L_{2}(\Omega)\rightarrow L_{2}(\Omega)}\leq\label{7.26}\\
&  \leq\|g_{2}\|_{L_{\infty}}\varepsilon^{\frac{1}{2}}\|[D_{2}%
,F(k)]F(k)\|_{L_{2}(\Omega)\rightarrow L_{2}(\Omega)}\|(S(k,\varepsilon
)+\varepsilon^{2}Q^{0})^{-1}\|_{L_{2}(\Omega)\rightarrow L_{2}(\Omega)}%
\leq\nonumber\\
&  \leq c_{1}\varepsilon^{\frac{1}{2}}C_{7}t(c_{0}k^{2}+c_{3}\varepsilon
^{2})^{-1}\leq c_{1}C_{7}c_{0}^{-\frac{1}{2}}c_{3}^{-\frac{1}{2}}%
\varepsilon^{-\frac{1}{2}}.\nonumber
\end{align}
Relations (\ref{1.1}), (\ref{6.13}) and (\ref{5.21}) imply that%
\begin{align}
\|\widehat{G}_{03}  &  (S(k,\varepsilon)P+\varepsilon^{2}Q^{0})^{-1}%
P\|_{L_{2}(\Omega)\rightarrow L_{2}(\Omega)}\leq\label{7.27}\\
&  \leq\|g_{2}\|_{L_{\infty}}\varepsilon^{\frac{1}{2}}\|[D_{2},F(k)]\|_{L_{2}%
(\Omega)\rightarrow L_{2}(\Omega)}\|F(k)^{\perp}(S(k,\varepsilon
)P+\varepsilon^{2}Q^{0})^{-1}P\|\leq\nonumber\\
&  \leq c_{1}\varepsilon^{\frac{1}{2}}C_{6}C_{4}c_{3}^{-\frac{1}{2}%
}\varepsilon^{-1}=c_{1}C_{6}C_{4}c_{3}^{-\frac{1}{2}}\varepsilon^{-\frac{1}%
{2}}.\nonumber
\end{align}
Finally, from (\ref{1.1}), (\ref{3.32}), (\ref{5.14}) and (\ref{5.3}) it
follows that%
\begin{align}
\|\overset{\vee}{G}_{03}  &  (S(k,\varepsilon)P+\varepsilon^{2}Q^{0}%
)^{-1}P\|_{L_{2}(\Omega)\rightarrow L_{2}(\Omega)}\leq\label{7.28}\\
&  \leq\|g_{2}\|_{L_{\infty}}\varepsilon^{-\frac{1}{2}}t\|\Phi(k)\|\cdot
\|\varepsilon D_{2}(S(k,\varepsilon)+\varepsilon^{2}Q^{0})^{-1}P\|_{L_{2}%
(\Omega)\rightarrow L_{2}(\Omega)}\leq\nonumber\\
&  \leq c_{1}\varepsilon^{-\frac{1}{2}}tC_{1}c_{0}^{-\frac{1}{2}}(c_{0}%
k^{2}+c_{3}\varepsilon^{2})^{-\frac{1}{2}}\leq c_{1}C_{1}c_{0}^{-1}%
\varepsilon^{-\frac{1}{2}}.\nonumber
\end{align}
As a result, relations (\ref{7.25})--(\ref{7.28}) imply that%
\begin{align}
\|G_{03}(S(k,\varepsilon)P+\varepsilon^{2}Q^{0})^{-1}P\|  &  _{L_{2}%
(\Omega)\rightarrow L_{2}(\Omega)}\leq C_{8}\varepsilon^{-\frac{1}{2}%
},\label{7.29}\\
C_{8}  &  =c_{1}(C_{7}c_{0}^{-\frac{1}{2}}c_{3}^{-\frac{1}{2}}+C_{6}C_{4}%
c_{3}^{-\frac{1}{2}}+C_{1}c_{0}^{-1}).\nonumber
\end{align}
Next, by (\ref{5.1}) and (\ref{5.12}), we have%
\begin{align}
\|G_{3}(A(k,\varepsilon)+\varepsilon^{2}Q)^{-1}\|_{L_{2}(\Omega)\rightarrow
L_{2}(\Omega)}  &  \leq\varepsilon^{\frac{1}{2}}\|\varepsilon D_{2}%
(A(k,\varepsilon)+\varepsilon^{2}Q)^{-1}\|_{L_{2}(\Omega)\rightarrow
L_{2}(\Omega)}\leq\label{7.30}\\
&  \leq c_{0}^{-\frac{1}{2}}\varepsilon^{\frac{1}{2}}(c_{0}k^{2}%
+c_{3}\varepsilon^{2})^{-\frac{1}{2}}\leq c_{0}^{-\frac{1}{2}}c_{3}^{-\frac
{1}{2}}\varepsilon^{-\frac{1}{2}}.\nonumber
\end{align}

\subsection{The form $\mathcal{I}_{4}$}

It remains to study the form (\ref{7.23}). It can be rewritten as%
\begin{align}
\mathcal{I}_{4}[u,v]  &  =(g_{2}^{0}\varepsilon D_{2}Pu,\varepsilon
D_{2}Pv)_{L_{2}(\Omega)}-(g_{2}\varepsilon D_{2}Pu,P\varepsilon D_{2}%
v)_{L_{2}(\Omega)}-\label{7.31}\\
&  -(g_{2}\varepsilon D_{2}Pu,P^{\perp}\varepsilon D_{2}v)_{L_{2}(\Omega
)}.\nonumber
\end{align}
Since $Pg_{2}P=g_{2}^{0}P$, the first two terms in the right-hand side of
(\ref{7.31}) compensate each other. Then%
\begin{equation}
\mathcal{I}_{4}[u,v]=-(g_{2}\varepsilon D_{2}Pu,P^{\perp}\varepsilon
D_{2}v)_{L_{2}(\mathbb{R}^{2})}=(G_{04}u,G_{4}v)_{L_{2}(\Omega)}, \label{7.32}%
\end{equation}
where $G_{04}=-g_{2}\varepsilon^{\frac{3}{2}}D_{2}P,\,G_{4}=\varepsilon
^{\frac{1}{2}}P^{\perp}D_{2}.$

By (\ref{1.1}), (\ref{5.3}) and (\ref{5.14}), we have
\begin{align}
\|G_{04}(S(k,\varepsilon)P+\varepsilon^{2}Q_{0})^{-1}P\|_{L_{2}(\Omega
)\rightarrow L_{2}(\Omega)}  &  \leq\|g_{2}\|_{L_{\infty}}\varepsilon
^{\frac{1}{2}}\|\varepsilon D_{2}(S(k,\varepsilon)P+\varepsilon^{2}Q_{0}%
)^{-1}P\|_{L_{2}(\Omega)\rightarrow L_{2}(\Omega)}\leq\label{7.33}\\
&  \leq c_{1}\varepsilon^{\frac{1}{2}}c_{0}^{-\frac{1}{2}}(c_{0}k^{2}%
+c_{3}\varepsilon^{2})^{-\frac{1}{2}}\leq c_{1}c_{0}^{-\frac{1}{2}}%
c_{3}^{-\frac{1}{2}}\varepsilon^{-\frac{1}{2}}.\nonumber
\end{align}
Next,%
\begin{align*}
G_{4}(A(k,\varepsilon)+\varepsilon^{2}Q)^{-1}  &  =\varepsilon^{\frac{1}{2}%
}P^{\perp}D_{2}(A(k,\varepsilon)+\varepsilon^{2}Q)^{-1}=\\
&  =\varepsilon^{\frac{1}{2}}P^{\perp}[D_{2},(A(k,\varepsilon)+\varepsilon
^{2}Q)^{-1}]+\varepsilon^{-\frac{1}{2}}P^{\perp}(A(k,\varepsilon
)+\varepsilon^{2}Q)^{-1}\varepsilon D_{2}.
\end{align*}
The first term on the right can be estimated with the help of (\ref{6.5}), and
the second one is estimated by using (\ref{5.12}) and (\ref{5.17}). As a
result, we obtain that
\begin{align}
\|G_{4}(A(k,\varepsilon)+\varepsilon^{2}Q)^{-1}\|_{L_{2}(\Omega)\rightarrow
L_{2}(\Omega)}  &  \leq\varepsilon^{\frac{1}{2}}C_{5}(c_{0}k^{2}%
+c_{3}\varepsilon^{2})^{-\frac{1}{2}} +\varepsilon^{-\frac{1}{2}}C_{3}%
c_{0}^{-\frac{1}{2}}\leq\label{7.34}\\
&  \leq(C_{5}c_{3}^{-\frac{1}{2}} + C_{3} c_{0}^{-\frac{1}{2}})\varepsilon
^{-\frac{1}{2}}.\nonumber
\end{align}

\subsection{Representation of the form $\mathcal{I}$}

Now we summarize the study of the form (\ref{7.4}). By (\ref{7.4}),
(\ref{7.8}), (\ref{7.12}), (\ref{7.15}), (\ref{7.18}), (\ref{7.21}),
(\ref{7.24}), (\ref{7.32}), this form is represented as%
\begin{align*}
\mathcal{I}[u,v]  &  =(G_{01}^{\prime}u,G_{1}^{\prime}v)_{L_{2}(\Omega
)}+(G_{01}^{\prime\prime}u,G_{1}^{\prime\prime}v)_{L_{2}(\Omega)}+\\
&  +(G_{02}u,G_{2}v)_{L_{2}(\Omega)}+(G_{03}u,G_{3}v)_{L_{2}(\Omega)}%
+(G_{04}u,G_{4}v)_{L_{2}(\Omega)}.
\end{align*}
We put $\mathfrak{G}=L_{2}(\Omega;\mathbb{C}^{5})=(L_{2}(\Omega))^{5}$ and
introduce the operators%
\begin{align*}
G_{0}  &  :L_{2}(\Omega)\rightarrow\mathfrak{G},\,\ G_{0}=\hbox {col}\{{G_{01}%
^{\prime},G_{01}^{\prime\prime},G_{02},G_{03},G_{04}\}},\\
G  &  :L_{2}(\Omega)\rightarrow\mathfrak{G},\,\ G=\hbox {col}\{{G_{1}^{\prime
},G_{1}^{\prime\prime},G_{2},G_{3},G_{4}\}},\\
\operatorname*{Dom}G_{0}  &  =\operatorname*{Dom}G=\{{u\in L_{2}(\Omega
):D_{2}u\in L_{2}(\Omega)\}}.
\end{align*}
The operators $G_{0}$ and $G$ are densely defined, and we have%
\[
\mathfrak{d}=\widetilde{H}^{1}(\Omega)\subset\operatorname*{Dom}%
G_{0}=\operatorname*{Dom}G,
\]
and%
\begin{equation}
\mathcal{I}[u,v]=(G_{0}u,Gv)_{\mathfrak{G}},\,\ u,v\in\mathfrak{d}.
\label{7.35}%
\end{equation}
Relations (\ref{7.14}), (\ref{7.17}), (\ref{7.20}), (\ref{7.30}), (\ref{7.34})
show that the operator $G(A(k,\varepsilon)+\varepsilon^{2}Q)^{-1}$ is bounded
and%
\begin{equation}
\|G(A(k,\varepsilon)+\varepsilon^{2}Q)^{-1}\|_{L_{2}(\Omega)\rightarrow
\mathfrak{G}}\leq C_{9}\varepsilon^{-\frac{1}{2}}, \label{7.36}%
\end{equation}
where%
\begin{equation}
C_{9}^{2}=C_{3}^{2}c_{3}^{-1}+c_{0}^{-\frac{1}{2}}c_{3}^{-\frac{3}{2}}%
+c_{0}^{-\frac{3}{2}}c_{3}^{-\frac{1}{2}}+c_{0}^{-1}c_{3}^{-1}+ (C_{5}%
c_{3}^{-\frac{1}{2}}+C_{3}c_{0}^{-\frac{1}{2}})^{2}. \label{7.37}%
\end{equation}
Inequalities (\ref{7.13}), (\ref{7.16}), (\ref{7.19}), (\ref{7.29}) and
(\ref{7.33}) imply that the operator $G_{0}(S(k,\varepsilon)P+\varepsilon
^{2}Q^{0})^{-1}P$ is bounded and%
\begin{equation}
\|G_{0}(S(k,\varepsilon)P+\varepsilon^{2}Q^{0})^{-1}P\|_{L_{2}(\Omega
)\rightarrow\mathfrak{G}}\leq C_{10}\varepsilon^{-\frac{1}{2}}, \label{7.38}%
\end{equation}
where%
\begin{equation}
C_{10}^{2}=c_{4}^{2}c_{3}^{{-2}}+C_{1}^{2}c_{4}^{2}c_{0}^{-\frac{1}{2}}%
c_{3}^{-\frac{3}{2}}+C_{2}^{2}c_{0}^{-\frac{3}{2}}c_{3}^{-\frac{1}{2}}%
+C_{8}^{2}+c_{1}^{2}c_{0}^{-1}c_{3}^{-1}. \label{7.39}%
\end{equation}
Thus, we have checked that in our case condition $5^{\circ}$ of Lemma
\ref{lem7.1} is satisfied.

Applying Lemma \ref{lem7.1}, we obtain representation of the form (\ref{7.3})
for the operator (\ref{5.23}):%
\begin{align}
{\mathcal{J}}(k,\varepsilon)  &  =(A(k,\varepsilon)+\varepsilon^{2}%
Q)^{-1}P-F(k)(S(k,\varepsilon)+\varepsilon^{2}Q^{0})^{-1}P=\label{7.40}\\
&  =(G(A(k,\varepsilon)+\varepsilon^{2}Q)^{-1})^{\ast}(G_{0}(S(k,\varepsilon
)P+\varepsilon^{2}Q^{0})^{-1}P).\nonumber
\end{align}
Then relations (\ref{7.36}), (\ref{7.38}) and (\ref{7.40}) imply that%
\[
\|{\mathcal{J}}(k,\varepsilon)\|_{L_{2}(\Omega)\rightarrow L_{2}(\Omega)}\leq
C_{9}C_{10}\varepsilon^{-1},\,\ \ |k|\leq t_{0},\ 0<\varepsilon\leq1.
\]
This proves estimate (\ref{5.24}) with $\widehat{C}=C_{9}C_{10}.$

\textit{This completes the proof of Proposition \textrm{\ref{prop5.5}}, and,
with it, of Theorem \textrm{\ref{theor4.1}}}. Also, \textit{this completes the
proof of Theorem} \ref{theor2.2} (see Subsection \ref{sect4.2}) \textit{and},
thereby (see Subsection \ref{sect2.2}), \textit{that of Theorem}
\ref{theor1.1}.

\section{Application to the Schr\"odinger operator with singular periodic
potential}

\label{sect10}

Homogenization problem for the Schr\"odinger operator $\mathcal{H}%
_{\varepsilon}= - \hbox{div} g(\mathbf{x}/\varepsilon) \nabla+ \varepsilon
^{-2} V(\mathbf{x}/\varepsilon)$ in $L_{2}(\mathbb{R}^{d})$ with coefficients
periodic in all directions was studied by a spectral method in Ch. 6, \S 1 of
\cite{BSu1}, \S 11 of \cite{BSu3}, \S 18 of \cite{BSu4}. Investigation was
based on a special factorization for the Schr\"odinger operator (see
\cite{KiSi} and \S 4 of \cite{BSu5}).

Here we apply Theorem 1 to homogenization problem for the Schr\"odinger
operator of the form
\begin{equation}
{H}_{\varepsilon}= D_{1} \widetilde{g}_{1}(\varepsilon^{-1}x_{1},x_{2})D_{1} +
D_{2} \widetilde{g}_{2}(\varepsilon^{-1}x_{1},x_{2})D_{2} + \varepsilon^{-2}
V(\varepsilon^{-1}x_{1},x_{2}) \label{S.0}%
\end{equation}
with metric and potential periodic in $x_{1}$. The operator ${H}_{\varepsilon
}$ contains a large factor $\varepsilon^{-2}$ standing at $V(\varepsilon
^{-1}x_{1},x_{2})$; in this sense the potential is singular. Under some
restrictions on the coefficients (potential $V$ must admit a certain
representation), we obtain approximation for the resolvent $({H}_{\varepsilon
}+ \lambda I)^{-1}$ in the operator $L_{2}$-norm with sharp order error estimate.

\subsection{Conditions on coefficients}

Let $\widetilde{g}_{1}(\mathbf{x})$, $\widetilde{g}_{2}(\mathbf{x})$ be
measurable functions in $\mathbb{R}^{2}$ satisfying the same conditions as
$g_{1}$, $g_{2}$. Namely, let $\widetilde{g}_{j}$, $j=1,2$, be periodic in
$x_{1}$ with period 1, uniformly bounded and positive definite:
\begin{equation}
\widetilde{c}_{0} \le\widetilde{g}_{j}(\mathbf{x}) \le\widetilde{c}_{1} <
\infty, \quad j=1,2,\ \hbox{a.~e.}\; \mathbf{x} \in\mathbb{R}^{2}. \label{S.1}%
\end{equation}
Also, it is assumed that $\widetilde{g}_{j}$ are Lipschitz functions with
respect to $x_{2}$:
\begin{equation}
\underset{\mathbf{x}\in{\mathbb{R}}^{2}}{\operatorname*{ess}\sup}|\partial_{2}
\widetilde{g}_{j}(\mathbf{x})|\leq\widetilde{c}_{2}<\infty,\quad j=1,2.
\label{S.2}%
\end{equation}

Next, let $\omega(\mathbf{x})$ be a measurable function in $\mathbb{R}^{2}$
such that
\begin{equation}
0< \omega_{0} \le\omega(\mathbf{x}) \le\omega_{1} < \infty,\quad
\hbox{a.~e.}\ \mathbf{x} \in\mathbb{R}^{2}. \label{S.3}%
\end{equation}
Assume that $\omega(\mathbf{x})$ is 1-periodic in $x_{1}$ and Lipschitz class
in $x_{2}$:
\begin{equation}
\underset{\mathbf{x}\in{\mathbb{R}}^{2}}{\operatorname*{ess}\sup}|\partial_{2}
\omega(\mathbf{x})|\leq\widetilde{c}_{3}<\infty. \label{S.4}%
\end{equation}

We put
\begin{equation}
V(\mathbf{x}) = - \frac{D_{1} \widetilde{g}_{1}(\mathbf{x})D_{1}
\omega(\mathbf{x})}{\omega(\mathbf{x})}, \label{S.5}%
\end{equation}
\begin{equation}
V_{2}(\mathbf{x}) = - \frac{D_{2} \widetilde{g}_{2}(\mathbf{x})D_{2}
\omega(\mathbf{x})}{\omega(\mathbf{x})}. \label{S.6}%
\end{equation}
Additional assumptions on $\widetilde{g}_{1}$, $\widetilde{g}_{2}$ and
$\omega$ are formulated in terms of the functions (\ref{S.5}), (\ref{S.6}).
Namely, assume that $V_{2}$ is uniformly bounded:
\begin{equation}
|V_{2}(\mathbf{x})| \le\widetilde{c}_{4},\quad\hbox{a.~e.}\ \mathbf{x}
\in\mathbb{R}^{2}, \label{S.7}%
\end{equation}
and Lipschitz class with respect to $x_{2}$:
\begin{equation}
\underset{\mathbf{x}\in{\mathbb{R}}^{2}}{\operatorname*{ess}\sup}|\partial
_{2}V_{2}(\mathbf{x})| \leq\widetilde{c}_{5}<\infty. \label{S.8}%
\end{equation}
Assume that $V(\cdot,x_{2}) \in L_{1}(0,1)$ for a.~e. $x_{2} \in\mathbb{R}$
and
\begin{equation}
\underset{x_{2}\in{\mathbb{R}}}{\operatorname*{ess}\sup}\|V(\cdot,
x_{2})\|_{L_{1}(0,1)} \leq\widetilde{c}_{6}<\infty. \label{S.9}%
\end{equation}

Using (\ref{S.1}), (\ref{S.3}), (\ref{S.5}) and (\ref{S.9}), it is easily seen
that $\omega(\mathbf{x})$ is Lipschitz in $x_{1}$:
\begin{equation}
\underset{\mathbf{x} \in{\mathbb{R}}^{2} }{\operatorname*{ess}\sup} |
\partial_{1} \omega(\mathbf{x})| \le\widetilde{c}_{7} = \omega_{1}
\widetilde{c}_{6}(\widetilde{c}_{0}^{-1} + \widetilde{c}_{1} \widetilde{c}%
_{0}^{-2})<\infty. \label{S.10}%
\end{equation}
Then, by (\ref{S.3}), (\ref{S.4}) and (\ref{S.10}), functions $\omega$ and
$\omega^{-1}$ are multipliers in the Sobolev class $H^{1}(\mathbb{R}^{2})$.

\subsection{Operators $\widetilde{A}_{\varepsilon}$, $A_{\varepsilon}$ and
$H_{\varepsilon}$}

For any measurable function $\varphi(x_{1},x_{2})$ which is 1-periodic in
$x_{1}$, we use the following notation
\[
\varphi^{\varepsilon}(x_{1},x_{2}) = \varphi(\varepsilon^{-1}x_{1},x_{2}).
\]

In $L_{2}(\mathbb{R}^{2})$, consider the operator $\widetilde{A}_{\varepsilon
}$ formally given by the expression
\begin{equation}
\widetilde{A}_{\varepsilon}= (\omega^{\varepsilon})^{-1}D_{1} \widetilde
{g}_{1}^{\varepsilon}(\omega^{\varepsilon})^{2}D_{1} (\omega^{\varepsilon
})^{-1} + (\omega^{\varepsilon})^{-1}D_{2} \widetilde{g}_{2}^{\varepsilon
}(\omega^{\varepsilon})^{2}D_{2} (\omega^{\varepsilon})^{-1}. \label{S.11}%
\end{equation}
Obviously,
\begin{equation}
\widetilde{A}_{\varepsilon}= (\omega^{\varepsilon})^{-1} A_{\varepsilon
}(\omega^{\varepsilon})^{-1}, \label{S.12}%
\end{equation}
where $A_{\varepsilon}$ is given by
\begin{equation}
A_{\varepsilon}= D_{1} g_{1}^{\varepsilon}D_{1} + D_{2} g_{2}^{\varepsilon
}D_{2} \label{S.12a}%
\end{equation}
(for precise definition of $A_{\varepsilon}$, see Section \ref{sect1}) with
the coefficients
\begin{equation}
g_{j}(\mathbf{x}) = \widetilde{g}_{j}(\mathbf{x}) \omega^{2}(\mathbf{x}),\quad
j=1,2. \label{S.13}%
\end{equation}
Due to conditions (\ref{S.1})--(\ref{S.4}) and periodicity of $\widetilde
{g}_{j}$, $\omega$ in $x_{1}$, coefficients (\ref{S.13}) satisfy the required
conditions (\ref{1.1})--(\ref{1.3}).

The precise definition of $\widetilde{A}_{\varepsilon}$ is given in terms of
the quadratic form
\begin{equation}
\widetilde{a}_{\varepsilon}[u,u] = \int_{\mathbb{R}^{2}} \left(
{g}^{\varepsilon}_{1} |D_{1} (\omega^{\varepsilon})^{-1}u|^{2} +{g}%
^{\varepsilon}_{2} |D_{2} (\omega^{\varepsilon})^{-1}u|^{2}\right)
\,d\mathbf{x},\quad u \in H^{1}(\mathbb{R}^{2}). \label{S.14}%
\end{equation}

By direct calculations, it is easy to check the following identity:
\begin{equation}
\widetilde{a}_{\varepsilon}[u,u] = \int_{\mathbb{R}^{2}} \left(  \widetilde
{g}_{1}^{\varepsilon}|D_{1} u|^{2} + \widetilde{g}_{2}^{\varepsilon}|D_{2}
u|^{2} + \varepsilon^{-2} V^{\varepsilon}|u|^{2} + V_{2}^{\varepsilon}%
|u|^{2}\right)  \,d\mathbf{x}, \quad u \in H^{1}(\mathbb{R}^{2}). \label{S.15}%
\end{equation}

Now, we introduce the Schr\"odinger operator ${H}_{\varepsilon}$ formally
given by (\ref{S.0}). The precise definition of ${H}_{\varepsilon}$ is given
in terms of the quadratic form
\begin{equation}
h_{\varepsilon}[u,u] = \int_{\mathbb{R}^{2}} \left(  \widetilde{g}%
_{1}^{\varepsilon}|D_{1} u|^{2} + \widetilde{g}_{2}^{\varepsilon}|D_{2} u|^{2}
+ \varepsilon^{-2} V^{\varepsilon}|u|^{2} \right)  \,d\mathbf{x}, \quad u \in
H^{1}(\mathbb{R}^{2}). \label{S.16}%
\end{equation}
Note that conditions (\ref{S.1}) and (\ref{S.9}) together with periodicity of
coefficients in $x_{1}$ ensure that the form ({\ref{S.16}) is closed and lower
semibounded on domain $H^{1}(\mathbb{R}^{2})$ (the same is true for the form
in the right-hand side of (\ref{S.15})). }

Relations (\ref{S.15}) and (\ref{S.16}) mean that
\begin{equation}
{H}_{\varepsilon}= \widetilde{A}_{\varepsilon}- V_{2}^{\varepsilon}.
\label{S.17}%
\end{equation}

\subsection{Homogenization for the Schr\"odinger operator $H_{\varepsilon}$}

From (\ref{S.17}) and (\ref{S.12}) it follows that
\begin{equation}
{H}_{\varepsilon}+ \lambda I = \widetilde{A}_{\varepsilon}- V_{2}%
^{\varepsilon}+ \lambda I = (\omega^{\varepsilon})^{-1} (A_{\varepsilon}+
Q_{\lambda}^{\varepsilon}) (\omega^{\varepsilon})^{-1},\quad\lambda
\in\mathbb{R}, \label{S.18}%
\end{equation}
where
\begin{equation}
Q_{\lambda}(\mathbf{x}) = (\lambda- V_{2}(\mathbf{x})) \omega^{2}(\mathbf{x}).
\label{S.19}%
\end{equation}
By (\ref{S.3}) and (\ref{S.7}), we have $Q_{\lambda}(\mathbf{x}) \ge
\lambda\omega_{0}^{2} - \omega_{1}^{2} \widetilde{c}_{4}$ for a.~e.
$\mathbf{x} \in\mathbb{R}^{2}$. Assume that $\lambda$ is sufficiently large so
that
\begin{equation}
\lambda\omega_{0}^{2} - \omega_{1}^{2} \widetilde{c}_{4} =: \widetilde
{c}_{\lambda}>0. \label{S.20}%
\end{equation}
Due to (\ref{S.3}), (\ref{S.4}), (\ref{S.7}), (\ref{S.8}), (\ref{S.20}) and
periodicity conditions, the function (\ref{S.19}) satisfies the same
conditions as $Q$ (see (\ref{1.6a})--(\ref{1.8})). We put
\begin{equation}
Q_{\lambda}^{0}(x_{2}) = \int_{0}^{1} Q_{\lambda}(x_{1},x_{2})\,dx_{1}.
\label{S.21}%
\end{equation}
Let $A^{0}$ be the effective operator corresponding to $A_{\varepsilon}$,
i.~e.,
\begin{equation}
A^{0} = D_{1} g_{1}^{0}(x_{2})D_{1} +D_{2} g_{2}^{0}(x_{2})D_{2}, \label{S.22}%
\end{equation}
where the coefficients $g_{1}^{0}$, $g_{2}^{0}$ are defined according to
(\ref{1.11}), (\ref{1.12}) in terms of the coefficients (\ref{S.13}).

By Theorem 1, we have
\begin{equation}
\|(A_{\varepsilon}+ Q_{\lambda}^{\varepsilon})^{-1} - (A^{0} +Q^{0}_{\lambda
})^{-1}\|_{L_{2}(\mathbb{R}^{2})\to L_{2}(\mathbb{R}^{2})} \le\widetilde
{C}_{\lambda}\varepsilon. \label{S.23}%
\end{equation}
Here the constant $\widetilde{C}_{\lambda}$ depends only on $\widetilde{c}%
_{0}$, $\widetilde{c}_{1}$, $\widetilde{c}_{2}$, $\omega_{0}$, $\omega_{1}$,
$\widetilde{c}_{3}$, $\widetilde{c}_{4}$, $\widetilde{c}_{5}$ and $\lambda$.
Multiplying operators under the norm sign in (\ref{S.23}) by $\omega
^{\varepsilon}$ from both sides and using identity (\ref{S.18}), we arrive at
the following result.

\begin{theorem}
\label{theorS.1}Let $\widetilde{g}_{1}$, $\widetilde{g}_{2}$, $\omega$ be
measurable functions in $\mathbb{R}^{2}$, periodic in $x_{1}$ with period $1$
and satisfying conditions \textrm{(\ref{S.1})--(\ref{S.4})}. Suppose that the
functions $V$, $V_{2}$ defined by \textrm{(\ref{S.5}), (\ref{S.6})} satisfy
conditions \textrm{(\ref{S.7})--(\ref{S.9})}. Let $H_{\varepsilon}$ be the
operator in $L_{2}(\mathbb{R}^{2})$ corresponding to the quadratic form
\textrm{(\ref{S.16})}. Let $A_{\varepsilon}$ be the operator
\textrm{(\ref{S.12a})} with the coefficients \textrm{(\ref{S.13})}, and let
$A^{0}$ be the corresponding effective operator \textrm{(\ref{S.22})} with
coefficients defined according to \textrm{(\ref{1.11}), (\ref{1.12})}. Let
$Q_{\lambda}$ be defined by \textrm{(\ref{S.19})}, and let restriction
\textrm{(\ref{S.20})} be satisfied. Let $Q_{\lambda}^{0}$ be given by
\textrm{(\ref{S.21})}. Then we have
\begin{equation}
\|(H_{\varepsilon}+\lambda I)^{-1} - \omega^{\varepsilon}(A^{0} +Q^{0}%
_{\lambda})^{-1} \omega^{\varepsilon}\|_{L_{2}(\mathbb{R}^{2})\to
L_{2}(\mathbb{R}^{2})} \le{C}_{\lambda}\varepsilon, \quad0 < \varepsilon\le1.
\label{S.24}%
\end{equation}
Here the constant $C_{\lambda}= \omega_{1}^{2} \widetilde{C}_{\lambda}$
depends only on $\widetilde{c}_{0}$, $\widetilde{c}_{1}$, $\widetilde{c}_{2}$,
$\omega_{0}$, $\omega_{1}$, $\widetilde{c}_{3}$, $\widetilde{c}_{4}$,
$\widetilde{c}_{5}$ and $\lambda$.
\end{theorem}

Note that the approximate operator in (\ref{S.24}) contains rapidly
oscillating factors $\omega^{\varepsilon}$ from both sides of the operator
$(A^{0} +Q^{0}_{\lambda})^{-1}$, but the inverse is taken only for the
differential operator $A^{0} +Q^{0}_{\lambda}$ with constant coefficients.
Such kind of results is typical also for the Schr\"odinger operator with
coefficients periodic in all directions (see \cite{BSu1}, Ch. 6, \S 1).

\subsection{About representation (\ref{S.5})}

A natural question is the following. If the initial object is the
Schr\"{o}dinger operator (\ref{S.0}) with given coefficients $\widetilde
{g}_{1}$, $\widetilde{g}_{2}$, $V$, then is it possible to find $\omega$ such
that representation (\ref{S.5}) is true and conditions of Theorem 18 are satisfied?

For a.e. $x_{2}\in\mathbb{R}$, consider the Schr\"{o}dinger operator
\begin{equation}
H(x_{2})=D_{1}\widetilde{g}_{1}(x_{1},x_{2})D_{1}+V(x_{1},x_{2}), \label{S.25}%
\end{equation}
acting in $L_{2}(\mathbb{R})$ and depending on parameter $x_{2}$. We assume
now that $\widetilde{g}_{1}$ and $V$ are 1-periodic in $x_{1}$ and satisfy
(\ref{S.1}) and (\ref{S.9}). Assume also that the \textit{bottom of the
spectrum of the operator} (\ref{S.25}) \textit{coincides with point}
$\lambda=0$:
\begin{equation}
\inf\text{spec}\,H(x_{2})=0. \label{S.26}%
\end{equation}
If initially condition (\ref{S.26}) was not satisfied, it is possible to
ensure this condition replacing $V(x_{1},x_{2})$ by $V(x_{1},x_{2}%
)-\lambda(x_{2})$, where $\lambda(x_{2})=\inf\text{spec}\,H(x_{2})$.

If condition (\ref{S.26}) is satisfied, there exists a positive periodic (in
$x_{1}$) solution $\omega(x_{1},x_{2})$ of the equation
\begin{equation}
D_{1} \widetilde{g}_{1}(x_{1},x_{2}) D_{1} \omega(x_{1},x_{2}) + V(x_{1}%
,x_{2}) \omega(x_{1},x_{2}) =0. \label{S.27}%
\end{equation}
This solution can be fixed by the condition
\begin{equation}
\int_{0}^{1} \omega^{2}(x_{1},x_{2}) \,dx_{1} =1. \label{S.28}%
\end{equation}
Moreover, this solution satisfies (\ref{S.3}), where $\omega_{0}$, $\omega
_{1}$ are controlled in terms of $\widetilde{c}_{0}$, $\widetilde{c}_{1}$,
$\widetilde{c}_{6}$. Concerning the properties of solution $\omega$, see
\cite{KiSi} and \S 4 of \cite{BSu5}. Thus, by (\ref{S.27}), $V$ satisfies
representation (\ref{S.5}).

Next, it is possible to ensure conditions (\ref{S.4}), (\ref{S.7}) and
(\ref{S.8}) imposing some smoothness assumptions on coefficients
$\widetilde{g}_{j}$, $V$ with respect to $x_{2}$. It suffices to assume that
derivatives $\partial_{2}^{l} \widetilde{g}_{1}$, $l=1,2,3,$ are uniformly
bounded, derivatives $\partial_{2}^{l} \widetilde{g}_{2}$, $l=1,2,$ are
uniformly bounded and that the norms $\|\partial_{2}^{l} V(\cdot
,x_{2})\|_{L_{1}(0,1)}$, $l=1,2,3,$ are uniformly bounded.

Under the above assumptions, all conditions of Theorem 18 are satisfied.

\section{Generalization of the main results for the case of arbitrary
dimension\label{sect11}}

\subsection{Generalization of Theorem 1}

One can study an analog of the operator $A_{\varepsilon}$ in $L_{2}%
(\mathbb{R}^{d})$ with $d=d_{1}+d_{2} \ge2$, $d_{1}, d_{2} \ge1$. We use the
notation $\mathbf{x} = (\mathbf{x}^{\prime}, \mathbf{x}^{\prime\prime})$,
$\mathbf{x}^{\prime}\in\mathbb{R}^{d_{1}}$, $\mathbf{x}^{\prime\prime}%
\in\mathbb{R}^{d_{2}}$, $\mathbf{D}_{\mathbf{x}^{\prime}} = -i \nabla
_{\mathbf{x}^{\prime}}$, $\mathbf{D}^{*}_{\mathbf{x}^{\prime}} = -i
\hbox{div}_{\mathbf{x}^{\prime}}$, and similarly for $\mathbf{D}%
_{\mathbf{x}^{\prime\prime}}$, $\mathbf{D}^{*}_{\mathbf{x}^{\prime\prime}}$.
Let $g_{1}(\mathbf{x}) = g_{1}(\mathbf{x}^{\prime}, \mathbf{x}^{\prime\prime
})$ be a $(d_{1} \times d_{1})$-matrix-valued measurable function and
$g_{2}(\mathbf{x}) = g_{2}(\mathbf{x}^{\prime}, \mathbf{x}^{\prime\prime})$ be
a $(d_{2} \times d_{2})$-matrix-valued measurable function in $\mathbb{R}^{d}%
$. Assume that $g_{1}(\mathbf{x})$, $g_{2}(\mathbf{x})$ are bounded and
positive definite:
\begin{equation}
c_{0} {\mathbf{1}}_{d_{1}} \le g_{1}(\mathbf{x}) \le c_{1} {\mathbf{1}}%
_{d_{1}}, \ \ c_{0} {\mathbf{1}}_{d_{2}} \le g_{2}(\mathbf{x}) \le c_{1}
{\mathbf{1}}_{d_{2}}, \quad0< c_{0} \le c_{1} < \infty
,\ \hbox{a.~e.}\ \mathbf{x} \in\mathbb{R}^{d}. \label{10.1}%
\end{equation}
Assume also that $g_{1}(\mathbf{x})$, $g_{2}(\mathbf{x})$ are periodic in
$\mathbf{x}^{\prime}$ with respect to some lattice $\Gamma\subset
\mathbb{R}^{d_{1}}$, and that $g_{1}$, $g_{2}$ are Lipschitz class with
respect to $\mathbf{x}^{\prime\prime}$:
\begin{equation}
\underset{ \mathbf{x} \in{\mathbb{R}}^{d}}{\operatorname*{ess}\sup}
|{\mathbf{D}}_{\mathbf{x}^{\prime\prime}} g_{j}(\mathbf{x})|\leq c_{2}%
<\infty,\quad j=1,2. \label{10.2}%
\end{equation}

In $L_{2}(\mathbb{R}^{d})$, consider the operator $A_{\varepsilon}$ formally
given by the differential expression
\begin{equation}
A_{\varepsilon}= {\mathbf{D}}^{*}_{\mathbf{x}^{\prime}} g_{1}\left(
\varepsilon^{-1}{\mathbf{x}^{\prime}}, {\mathbf{x}}^{\prime\prime}\right)
{\mathbf{D}}_{\mathbf{x}^{\prime}} + {\mathbf{D}}^{*}_{\mathbf{x}%
^{\prime\prime}} g_{1}\left(  \varepsilon^{-1}{\mathbf{x}^{\prime}},
{\mathbf{x}^{\prime\prime}} \right)  {\mathbf{D}}_{\mathbf{x}^{\prime\prime}%
},\quad\varepsilon>0. \label{10.2a}%
\end{equation}
Precisely, $A_{\varepsilon}$ is a selfadjoint operator in $L_{2}%
(\mathbb{R}^{d})$ corresponding to the quadratic form
\begin{equation}
a_{\varepsilon}[u,u]=\int_{{\mathbb{R}}^{d}}\left(  \left\langle g_{1}\left(
\varepsilon^{-1}{\mathbf{x}^{\prime}}, \mathbf{x}^{\prime\prime}\right)
{\mathbf{D}}_{\mathbf{x}^{\prime}} u,{\mathbf{D}}_{\mathbf{x}^{\prime}} u
\right\rangle _{{\mathbb{C}}^{d_{1}}} + \left\langle g_{2}\left(
\varepsilon^{-1}{\mathbf{x}^{\prime}}, \mathbf{x}^{\prime\prime}\right)
{\mathbf{D}}_{\mathbf{x}^{\prime\prime}} u,{\mathbf{D}}_{\mathbf{x}%
^{\prime\prime}} u \right\rangle _{{\mathbb{C}}^{d_{2}}} \right)
\,d\mathbf{x}, \quad u \in H^{1}({\mathbb{R}}^{d}). \label{10.3}%
\end{equation}

Let $Q(\mathbf{x}) = Q(\mathbf{x}^{\prime},\mathbf{x}^{\prime\prime})$ be a
real-valued measurable function in $\mathbb{R}^{d}$ such that
\begin{equation}
0< c_{3} \le Q(\mathbf{x}) \le c_{4} < \infty,\ \ \hbox{a.~e.}\ \mathbf{x}
\in\mathbb{R}^{d}. \label{10.4}%
\end{equation}
Assume that $Q(\mathbf{x})$ is $\Gamma$-periodic in $\mathbf{x}^{\prime}$, and
that
\begin{equation}
\underset{ \mathbf{x} \in{\mathbb{R}}^{d}}{\operatorname*{ess}\sup}
|{\mathbf{D}}_{\mathbf{x}^{\prime\prime}} Q(\mathbf{x})| \leq c_{5}<\infty.
\label{10.5}%
\end{equation}
We put
\begin{equation}
Q^{\varepsilon}(\mathbf{x}^{\prime}, \mathbf{x}^{\prime\prime}) = Q\left(
\varepsilon^{-1} {\mathbf{x}^{\prime}}, \mathbf{x}^{\prime\prime}\right)  .
\label{10.6}%
\end{equation}

For the generalized resolvent $(A_{\varepsilon}+ Q^{\varepsilon})^{-1}$, the
analog of Theorem \ref{theor1.1} is true. The effective operator $A^{0}$
corresponds to the differential expression
\begin{equation}
A^{0} = {\mathbf{D}}^{*}_{{\mathbf{x}}^{\prime}} g_{1}^{0}(\mathbf{x}%
^{\prime\prime}){\mathbf{D}}_{{\mathbf{x}}^{\prime}} + {\mathbf{D}}%
^{*}_{{\mathbf{x}}^{\prime\prime}} g_{2}^{0}(\mathbf{x}^{\prime\prime
}){\mathbf{D}}_{{\mathbf{x}}^{\prime\prime}}. \label{10.7}%
\end{equation}
Here the effective coefficients $g_{1}^{0}({\mathbf{x}}^{\prime\prime})$,
$g_{2}^{0}({\mathbf{x}}^{\prime\prime})$ are defined as follows. Let
$\Omega^{\prime}$ be an elementary cell of the lattice $\Gamma\subset
{\mathbb{R}}^{d_{1}}$. Then
\begin{equation}
g_{2}^{0}({\mathbf{x}}^{\prime\prime}) = |\Omega^{\prime-1} \int
_{\Omega^{\prime}} g_{2}(\mathbf{x}^{\prime},\mathbf{x}^{\prime\prime
})\,d\mathbf{x}^{\prime}. \label{10.8}%
\end{equation}
The matrix $g_{1}^{0}(\mathbf{x}^{\prime\prime})$ is the effective matrix
corresponding to the elliptic operator ${\mathbf{D}}^{*}_{\mathbf{x}^{\prime}}
g_{1}(\mathbf{x}^{\prime}, \mathbf{x}^{\prime\prime}){\mathbf{D}}%
_{\mathbf{x}^{\prime}}$ in $L_{2}(\mathbb{R}^{d_{1}})$. Recall the definition
of $g_{1}^{0}$. Let $\mathbf{e}_{1},\dots, \mathbf{e}_{d_{1}}$ be the standard
orthonormal basis in ${\mathbb{C}}^{d_{1}}$, and let $\Phi_{j}(\mathbf{x})$ be
a $\Gamma$-periodic (in $\mathbf{x}^{\prime}$) solution of the equation
\[
\hbox{div}_{\mathbf{x}^{\prime}} g_{1}(\mathbf{x}^{\prime},\mathbf{x}%
^{\prime\prime}) (\nabla_{\mathbf{x}^{\prime}} \Phi_{j}(\mathbf{x}^{\prime
},\mathbf{x}^{\prime\prime}) + \mathbf{e}_{j}) =0,
\]
$j=1,\dots, d_{1}$. Here $\mathbf{x}^{\prime\prime}$ is a parameter. Then
$g^{0}_{1}(\mathbf{x}^{\prime\prime})$ is a $(d_{1} \times d_{1})$-matrix with
the columns
\begin{equation}
{\mathbf{g}}_{j}(\mathbf{x}^{\prime\prime}) = |\Omega^{\prime-1} \int
_{\Omega^{\prime}} g_{1}(\mathbf{x}^{\prime},\mathbf{x}^{\prime\prime})
(\nabla_{\mathbf{x}^{\prime}} \Phi_{j}(\mathbf{x}^{\prime},\mathbf{x}%
^{\prime\prime}) + \mathbf{e}_{j})\,d \mathbf{x}^{\prime},\quad j=1,\dots,
d_{1}. \label{10.9}%
\end{equation}
Note that for $d_{1}=1$ relation (\ref{10.9}) is equivalent to (\ref{1.11}).

Finally, let
\begin{equation}
Q^{0}(\mathbf{x}^{\prime\prime}) = |\Omega^{\prime-1} \int_{\Omega^{\prime}}
Q(\mathbf{x}^{\prime}, \mathbf{x}^{\prime\prime})\,d\mathbf{x}^{\prime}.
\label{10.10}%
\end{equation}

The following result can be proved by analogy with the proof of Theorem
\ref{theor1.1}.

\begin{theorem}
\label{theor10.1}Let $d = d_{1} + d_{2} \ge2$, $d_{1}, d_{2} \ge1$. Let
$g_{1}(\mathbf{x})$ and $g_{2}(\mathbf{x})$ be measurable matrix-valued
functions in $\mathbb{R}^{d}$ \textrm{(}of sizes $d_{1}\times d_{1}$ and
$d_{2} \times d_{2}$ correspondingly\textrm{)}. Suppose that $g_{1}, g_{2}$
are $\Gamma$-periodic in $\mathbf{x}^{\prime}$ and satisfy conditions
\textrm{(\ref{10.1}), (\ref{10.2})}. Let $A_{\varepsilon}$ be the operator in
$L_{2}({\mathbb{R}}^{d})$ corresponding to the form \textrm{(\ref{10.3})}. Let
$Q(\mathbf{x})$ be a measurable function in $\mathbb{R}^{d}$, $\Gamma
$-periodic in $\mathbf{x}^{\prime}$ and satisfying conditions
\textrm{(\ref{10.4}), (\ref{10.5})}. Let $Q^{\varepsilon}(\mathbf{x})$ be
defined by \textrm{(\ref{10.6})}. Let $A^{0}$ be the operator
\textrm{(\ref{10.7})}, where the coefficients $g_{1}^{0}$, $g_{2}^{0}$ are
defined by \textrm{(\ref{10.9})} and \textrm{(\ref{10.8})} respectively. Let
$Q^{0}$ be defined by \textrm{(\ref{10.10})}. Then we have%
\begin{equation}
\|(A_{\varepsilon}+Q^{\varepsilon})^{-1}-(A^{0}+Q^{0})^{-1}\|_{L_{2}%
({\mathbb{R}}^{d})\rightarrow L_{2}({\mathbb{R}}^{d})}\leq C\varepsilon
,\ \ 0<\varepsilon\leq1, \label{10.11}%
\end{equation}
where the constant $C$ depends only on $c_{j},\,j=0,...,5,$ and on parameters
of the lattice $\Gamma$.
\end{theorem}

\subsection{Generalization of Theorem 18}

Let $\widetilde{g}_{j}(\mathbf{x})$, $j=1,2,$ be measurable matrix-valued
functions in $\mathbb{R}^{d}$ satisfying the same conditions as $g_{j}$.
Namely, we assume that $\widetilde{g}_{j}$, $j=1,2,$ are $\Gamma$-periodic in
$\mathbf{x}^{\prime}$, bounded and positive definite:
\begin{equation}
\widetilde{c}_{0} {\mathbf{1}}_{d_{1}} \le\widetilde{g}_{1}(\mathbf{x})
\le\widetilde{c}_{1} {\mathbf{1}}_{d_{1}}, \ \ \widetilde{c}_{0} {\mathbf{1}%
}_{d_{2}} \le\widetilde{g}_{2}(\mathbf{x}) \le\widetilde{c}_{1} {\mathbf{1}%
}_{d_{2}}, \quad0< c_{0} \le c_{1} < \infty,\ \hbox{a.~e.}\ \mathbf{x}
\in\mathbb{R}^{d}. \label{10.12}%
\end{equation}
Assume also that $\widetilde{g}_{1}(\mathbf{x})$, $\widetilde{g}%
_{2}(\mathbf{x})$ are Lipschitz with respect to $\mathbf{x}^{\prime\prime}$:
\begin{equation}
\underset{ \mathbf{x} \in{\mathbb{R}}^{d}}{\operatorname*{ess}\sup}
|{\mathbf{D}}_{\mathbf{x}^{\prime\prime}} \widetilde{g}_{j}(\mathbf{x})|\leq
c_{2}<\infty,\quad j=1,2. \label{10.13}%
\end{equation}

Next, let $\omega(\mathbf{x})$ be a measurable function in $\mathbb{R}^{d}$
such that
\begin{equation}
0< \omega_{0} \le\omega(\mathbf{x}) \le\omega_{1} < \infty,\quad
\hbox{a.~e.}\ \mathbf{x} \in\mathbb{R}^{d}. \label{10.14}%
\end{equation}
Assume that $\omega(\mathbf{x})$ is $\Gamma$-periodic in $\mathbf{x}^{\prime}$
and Lipschitz class in $\mathbf{x}^{\prime\prime}$:
\begin{equation}
\underset{\mathbf{x}\in{\mathbb{R}}^{2}}{\operatorname*{ess}\sup}%
|\mathbf{D}_{\mathbf{x}^{\prime\prime}} \omega(\mathbf{x})|\leq\widetilde
{c}_{3}<\infty. \label{10.15}%
\end{equation}

We put
\begin{equation}
V(\mathbf{x}) = - \frac{\mathbf{D}_{\mathbf{x}^{\prime}}^{*} \widetilde{g}%
_{1}(\mathbf{x})\mathbf{D}_{\mathbf{x}^{\prime}} \omega(\mathbf{x})}%
{\omega(\mathbf{x})}, \label{10.16}%
\end{equation}
\begin{equation}
V_{2}(\mathbf{x}) = - \frac{\mathbf{D}_{\mathbf{x}^{\prime\prime}}^{*}
\widetilde{g}_{2}(\mathbf{x}) \mathbf{D}_{\mathbf{x}^{\prime\prime}}
\omega(\mathbf{x})}{\omega(\mathbf{x})}. \label{10.17}%
\end{equation}
Additional assumptions on $\widetilde{g}_{1}$, $\widetilde{g}_{2}$ and
$\omega$ are formulated in terms of the functions (\ref{10.16}),
(\ref{10.17}). Namely, assume that $V_{2}$ is uniformly bounded:
\begin{equation}
|V_{2}(\mathbf{x})| \le\widetilde{c}_{4},\quad\hbox{a.~e.}\ \mathbf{x}
\in\mathbb{R}^{d}, \label{10.18}%
\end{equation}
and Lipschitz class with respect to $\mathbf{x}^{\prime\prime}$:
\begin{equation}
\underset{\mathbf{x}\in{\mathbb{R}}^{2}}{\operatorname*{ess}\sup}%
|\mathbf{D}_{\mathbf{x}^{\prime\prime}} V_{2}(\mathbf{x})| \leq\widetilde
{c}_{5}<\infty. \label{10.19}%
\end{equation}
Assume that $V(\cdot,\mathbf{x}^{\prime\prime}) \in L_{p}(\Omega^{\prime})$
for a.~e. $\mathbf{x}^{\prime\prime}\in\mathbb{R}^{d_{2}}$, where $p=1$ for
$d_{1}=1$ and $2p>d_{1}$ for $d_{1} \ge2$, and
\begin{equation}
\underset{\mathbf{x}^{\prime\prime}\in{\mathbb{R}^{d_{2}}}}%
{\operatorname*{ess}\sup}\|V(\cdot, \mathbf{x}^{\prime\prime}) \|_{L_{p}%
(\Omega^{\prime})} \leq\widetilde{c}_{6}<\infty. \label{10.20}%
\end{equation}

We put
\begin{equation}
g_{j}(\mathbf{x}) = \widetilde{g}_{j}(\mathbf{x}) \omega^{2}(\mathbf{x}),\quad
j=1,2, \label{10.21}%
\end{equation}
\begin{equation}
Q_{\lambda}(\mathbf{x}) = (\lambda- V_{2}(\mathbf{x})) \omega^{2}(\mathbf{x}),
\label{10.22}%
\end{equation}
\begin{equation}
Q_{\lambda}^{0}(\mathbf{x}^{\prime\prime}) = |\Omega^{\prime-1}\int
_{\Omega^{\prime}} Q_{\lambda}(\mathbf{x}^{\prime}, \mathbf{x}^{\prime\prime
})\,d\mathbf{x}^{\prime}. \label{10.23}%
\end{equation}
Suppose that $\lambda$ is such that
\begin{equation}
\lambda\omega^{2}_{0} - \omega_{1}^{2} \widetilde{c}_{4} >0. \label{10.24}%
\end{equation}

For any $\Gamma$-periodic (in $\mathbf{x}^{\prime}$) function $\varphi$, we
use the notation
\[
\varphi^{\varepsilon}(\mathbf{x}) = \varphi(\varepsilon^{-1} \mathbf{x}%
^{\prime}, \mathbf{x}^{\prime\prime}).
\]

In $L_{2}(\mathbb{R}^{d})$, consider the Schr\"odinger operator
\begin{equation}
H_{\varepsilon}= \mathbf{D}^{*}_{\mathbf{x}^{\prime}} \widetilde{g}%
_{1}^{\varepsilon}\mathbf{D}_{\mathbf{x}^{\prime}} + \mathbf{D}^{*}%
_{\mathbf{x}^{\prime\prime}} \widetilde{g}_{2}^{\varepsilon}\mathbf{D}%
_{\mathbf{x}^{\prime\prime}} + \varepsilon^{-2} V^{\varepsilon}. \label{10.25}%
\end{equation}
The precise definition of $H_{\varepsilon}$ is given in terms of the quadratic
form
\begin{equation}
h_{\varepsilon}[u,u]= \int_{\mathbb{R}^{d}} \left(  \langle\widetilde
{g}^{\varepsilon}_{1} \mathbf{D}_{\mathbf{x}^{\prime}}u, \mathbf{D}%
_{\mathbf{x}^{\prime}}u\rangle_{\mathbb{C}^{d_{1}}} + \langle\widetilde
{g}^{\varepsilon}_{2} \mathbf{D}_{\mathbf{x}^{\prime\prime}}u, \mathbf{D}%
_{\mathbf{x}^{\prime\prime}}u\rangle_{\mathbb{C}^{d_{2}}}+ \varepsilon
^{-2}V^{\varepsilon}|u|^{2}\right)  \,d\mathbf{x},\quad u \in H^{1}%
(\mathbb{R}^{d}). \label{10.26}%
\end{equation}

The following result can be deduced from Theorem 19 similarly to deduction of
Theorem 18 from Theorem 1.

\begin{theorem}
\label{theor10.2}Let $\widetilde{g}_{1}$, $\widetilde{g}_{2}$, $\omega$ be
measurable functions in $\mathbb{R}^{d}$, periodic in $\mathbf{x}^{\prime}$
with respect to the lattice $\Gamma$ and satisfying conditions
\textrm{(\ref{10.12})--(\ref{10.15})}. Suppose that the functions $V$, $V_{2}$
defined by \textrm{(\ref{10.16}), (\ref{10.17})} satisfy conditions
\textrm{(\ref{10.18})--(\ref{10.20})}. Let $H_{\varepsilon}$ be the operator
in $L_{2}(\mathbb{R}^{d})$ corresponding to the quadratic form
\textrm{(\ref{10.26})}. Let $A_{\varepsilon}$ be the operator
\textrm{(\ref{10.2a})} with the coefficients \textrm{(\ref{10.21})}, and let
$A^{0}$ be the corresponding effective operator \textrm{(\ref{10.7})} with
coefficients defined according to \textrm{(\ref{10.8}), (\ref{10.9})}. Let
$Q_{\lambda}$ be defined by \textrm{(\ref{10.22})}, and let restriction
\textrm{(\ref{10.24})} be satisfied. Let $Q_{\lambda}^{0}$ be given by
\textrm{(\ref{10.23})}. Then we have
\begin{equation}
\|(H_{\varepsilon}+\lambda I)^{-1} - \omega^{\varepsilon}(A^{0} +Q^{0}%
_{\lambda})^{-1} \omega^{\varepsilon}\|_{L_{2}(\mathbb{R}^{d})\to
L_{2}(\mathbb{R}^{d})} \le{C}_{\lambda}\varepsilon, \quad0 < \varepsilon\le1.
\label{10.27}%
\end{equation}
Here the constant $C_{\lambda}$ depends only on $\widetilde{c}_{0}$,
$\widetilde{c}_{1}$, $\widetilde{c}_{2}$, $\omega_{0}$, $\omega_{1}$,
$\widetilde{c}_{3}$, $\widetilde{c}_{4}$, $\widetilde{c}_{5}$, $\lambda$ and
parameters of the lattice $\Gamma$.
\end{theorem}

Similarly to Subsection 10.4, it is possible to give sufficient conditions on
$\widetilde{g}_{1}$, $\widetilde{g}_{2}$, $V$, which ensure representation
(\ref{10.16}) and conditions of Theorem 20. For this, consider the
Schr\"odinger operator
\begin{equation}
H(\mathbf{x}^{\prime\prime})= \mathbf{D}_{\mathbf{x}^{\prime}}^{*}
\widetilde{g}_{1}(\mathbf{x}^{\prime},\mathbf{x}^{\prime\prime})
\mathbf{D}_{\mathbf{x}^{\prime}} + \varepsilon^{-2} V(\mathbf{x}^{\prime},
\mathbf{x}^{\prime\prime}), \label{10.28}%
\end{equation}
acting in $L_{2}(\mathbb{R}^{d_{1}})$ and depending on parameter
$\mathbf{x}^{\prime\prime}\in\mathbb{R}^{d_{2}}$. Assume that $\widetilde
{g}_{j}$, $V$ are $\Gamma$-periodic in $\mathbf{x}^{\prime}$ and satisfy
(\ref{10.12}) and (\ref{10.20}). Assume also that the \textit{bottom of the
spectrum of the operator} (\ref{10.28}) \textit{coincides with point}
$\lambda=0$:
\begin{equation}
\inf\text{spec}\,H(\mathbf{x}^{\prime\prime}) =0,\quad\text{a.~e.}%
\ \mathbf{x}^{\prime\prime}\in\mathbb{R}^{d_{2}}. \label{10.29}%
\end{equation}

If condition (\ref{10.29}) is satisfied, there exists a positive periodic (in
$\mathbf{x}^{\prime}$) solution $\omega(\mathbf{x}^{\prime},\mathbf{x}%
^{\prime\prime})$ of the equation
\begin{equation}
\mathbf{D}_{\mathbf{x}^{\prime}}^{*} \widetilde{g}_{1}(\mathbf{x}^{\prime
},\mathbf{x}^{\prime\prime}) \mathbf{D}_{\mathbf{x}^{\prime}} \omega
(\mathbf{x}^{\prime},\mathbf{x}^{\prime\prime}) + V(\mathbf{x}^{\prime
},\mathbf{x}^{\prime\prime}) \omega(\mathbf{x}^{\prime},\mathbf{x}%
^{\prime\prime}) =0. \label{10.30}%
\end{equation}
This solution can be fixed by the condition
\begin{equation}
\int_{\Omega^{\prime}} \omega^{2}(\mathbf{x}^{\prime},\mathbf{x}^{\prime
\prime}) \,d\mathbf{x}^{\prime}=|\Omega^{\prime}|. \label{10.31}%
\end{equation}
Moreover, this solution satisfies (\ref{10.14}), where $\omega_{0}$,
$\omega_{1}$ are controlled in terms of $\widetilde{c}_{0}$, $\widetilde
{c}_{1}$, $\widetilde{c}_{6}$. Concerning the properties of solution $\omega$,
see \cite{KiSi} and \S 4 of \cite{BSu5}. By (\ref{10.30}), representation
(\ref{10.16}) is true.

Next, it is possible to ensure conditions (\ref{10.15}), (\ref{10.18}) and
(\ref{10.19}) imposing some smoothness assumptions on coefficients
$\widetilde{g}_{j}$, $V$ with respect to $\mathbf{x}^{\prime\prime}$. It
suffices to assume that derivatives $D_{\mathbf{x}^{\prime\prime}}^{\alpha
}\widetilde{g}_{1}$, $|\alpha|\le3,$ are uniformly bounded, derivatives
$D_{\mathbf{x}^{\prime\prime}}^{\alpha}\widetilde{g}_{2}$, $|\alpha|\le2$, are
uniformly bounded and that the norms $\|D_{\mathbf{x}^{\prime\prime}}^{\alpha
}V(\cdot,\mathbf{x}^{\prime\prime}) \|_{L_{p}(\Omega^{\prime})}$, $|\alpha
|\le3,$ are uniformly bounded.

Under the above assumptions, all conditions of Theorem 20 are satisfied.

\end{document}